\documentclass{amsart}
\usepackage{amssymb}
\usepackage{amscd}
\usepackage{verbatim}
\usepackage{epsfig}
\begin{document}

\newcommand\Mand{\ \text{and}\ }
\newcommand\Mwith{\ \text{with}\ }
\newcommand\Mfor{\ \text{for}\ }
\newcommand\Mst{\ \text{such that}\ }
\newcommand\Mor{\ \text{or}\ }
\newcommand\Mif{\ \text{if}\ }
\newcommand\Miff{\ \text{iff}\ }
\newcommand\Mthen{\ \text{then}\ }
\newcommand\nin{\notin}
\newcommand\identity{\operatorname{id}}
\newcommand\Id{\operatorname{Id}}
\newcommand\Real{\mathbb{R}}
\newcommand\pos{\Real^+}
\newcommand\Rnp{\Real\setminus\{0\}}
\newcommand\nzero{\setminus\{0\}}
\newcommand\Cx{\mathbb{C}}
\newcommand\Cxp{\Cx^+}
\newcommand\Cxm{\Cx^-}
\newcommand\Nat{\mathbb{N}}
\newcommand\halfNat{{\frac{1}{2}}\mathbb{N}}
\newcommand\intgr{\mathbb{Z}}
\newcommand\im{\operatorname{Im}}
\newcommand\re{\operatorname{Re}}
\newcommand\sign{\operatorname{sign}}
\newcommand\codim{\operatorname{codim}}
\newcommand\End{\operatorname{End}}
\newcommand\Ker{\operatorname{Ker}}
\newcommand\Hom{\operatorname{Hom}}
\newcommand\tr{\operatorname{tr}}
\newcommand\Tr{\operatorname{Tr}}
\newcommand\ideal{{\mathcal I}}
\newcommand\Span{\operatorname{span}}
\newcommand\image{\operatorname{image}}
\newcommand\Range{\operatorname{Ran}}
\newcommand\Graph{\operatorname{graph}}
\newcommand\slim{\operatornamewithlimits{s-lim}}
\newcommand\sll{\operatorname{sl}}
\newcommand\sol{\operatorname{so}}
\newcommand\SL{\operatorname{SL}}
\newcommand\SO{\operatorname{SO}}
\newcommand\On{\operatorname{O}}
\newcommand\pa{\partial}
\newcommand\Rn{\Real^n}
\newcommand\Rm{\Real^m}
\newcommand\RN{\Real^N}
\newcommand\RtN{\Real^{2N}}
\newcommand\RM{\Real^M}
\newcommand\sphere{\mathbb{S}}
\newcommand\Sn{\sphere^{n-1}}
\newcommand\Sm{\sphere^{m-1}}
\newcommand\Snp{\sphere^n_+}
\newcommand\Smp{\sphere^m_+}
\newcommand\SN{\sphere^{N-1}}
\newcommand\SNp{\sphere^N_+}
\newcommand\circlep{\sphere^1_+}
\newcommand\Phom{P_{h}}
\newcommand\Shom{S_{h}}
\newcommand\distance{\operatorname{dist}}
\newcommand\cl{\operatorname{cl}}
\newcommand\interior{\operatorname{int}}
\newcommand\Fa{\operatorname{Fa}}
\newcommand\ff{\operatorname{ff}}
\newcommand\mf{\operatorname{mf}}
\newcommand\cf{\operatorname{cf}}
\newcommand\scf{\operatorname{sf}}
\newcommand\lf{\operatorname{lf}}
\newcommand\rf{\operatorname{rf}}
\newcommand\indfam{{\mathcal K}}
\newcommand\fraka{{\mathfrak a}}
\newcommand\calA{{\mathcal A}}
\newcommand\calB{{\mathcal B}}
\newcommand\calR{{\mathcal R}}
\newcommand\calO{{\mathcal O}}
\newcommand\calJ{{\mathcal J}}
\newcommand\calM{{\mathcal M}}
\newcommand\calN{{\mathcal N}}
\newcommand\calX{{\mathcal X}}
\newcommand\calF{{\mathcal F}}
\newcommand\calG{{\mathcal G}}
\newcommand\calT{{\mathcal T}}
\newcommand\calC{{\mathcal C}}
\newcommand\calCt{{\tilde {\mathcal C}}}
\newcommand\calCL{{\mathcal C}_{\text L}}
\newcommand\calCR{{\mathcal C}_{\text R}}
\newcommand\Cinf{{\mathcal C}^{\infty}}
\newcommand\dist{{\mathcal C}^{-\infty}}
\newcommand\dCinf{\dot\Cinf}
\newcommand\ddist{\dot\dist}
\newcommand\Cj{{\mathcal C}^j}
\newcommand\Linf{L^{\infty}}
\newcommand\phg{{\text{phg}}}
\newcommand\bcon{{\mathcal A}}
\newcommand\bconc{{\mathcal A}_{\text{phg}}}
\newcommand\Sch{{\mathcal S}}
\newcommand\temp{\Sch^{\prime}}
\newcommand\Diff{\operatorname{Diff}}
\newcommand\Diffb{\operatorname{Diff}_{\text{b}}}
\newcommand\Diffc{\operatorname{Diff}_{\text{c}}}
\newcommand\Diffsc{\operatorname{Diff}_{\text{sc}}}
\newcommand\DiffI{\operatorname{Diff}_{\text{I}}}
\newcommand\DiffIq{\operatorname{Diff}_{\text{I},q}}
\newcommand\sing{\text{sing}}
\newcommand\reg{\text{reg}}
\newcommand\supp{\operatorname{supp}}
\newcommand\ssupp{\operatorname{sing\ supp}}
\newcommand\csupp{\operatorname{cone\ supp}}
\newcommand\esupp{\operatorname{ess\ supp}}
\newcommand\Fr{{\mathcal F}}
\newcommand\Frinv{\Fr^{-1}}
\newcommand\bop{{\mathcal B}}
\newcommand\spec{\operatorname{spec}}
\newcommand\pspec{\spec_{pp}}
\newcommand\cspec{\spec_{c}}
\newcommand\FIO{{\mathcal I}}
\newcommand\SP{\operatorname{RC}}
\newcommand\RC{\operatorname{RC}}
\newcommand\Symc{S_c}
\newcommand\Symca{S_c^{\alpha}}
\newcommand\Symczero{S_c^{0,...,0}}
\newcommand\sci{{}^{\text{sc}}}
\newcommand\sct{\sci T^*}
\newcommand\scdt{\sci \dot T^*}
\newcommand\dS{\dot S^*}
\newcommand\dT{\dot T^*}
\newcommand\dSreg{\dot\Sigma_{\text reg}}
\newcommand\scct{\sci\bar{T}^*}
\newcommand\Csc{C_{\text{sc}}}
\newcommand\SNpscd{(\SNp)^2_{\text{sc}}}
\newcommand\scdiag{\Delta_{\text{sc}}}
\newcommand\projscl{\pi^L_{\text{sc}}}
\newcommand\projscr{\pi^R_{\text{sc}}}
\newcommand\scHL{\sci H^{2,0}_{|\zeta|^2-\lambda^2}}
\newcommand\scHrg{\sci H^{2,0}_{\sqrt{g}}}
\newcommand\Hsc{H_{\text{sc}}}
\newcommand\WF{\operatorname{WF}}
\newcommand\WFp{\operatorname{WF^{\prime}}}
\newcommand\WFsc{\operatorname{WF}_{\text{sc}}}
\newcommand\WFscp{\operatorname{WF_{sc}^{\prime}}}
\newcommand\WFC{\operatorname{WF}_C}
\newcommand\WFCi{\operatorname{WF}_{C_i}}
\newcommand\elliptic{\operatorname{ell}}
\newcommand\Psop{\operatorname{\Psi}}
\newcommand\Psiscrs{\operatorname{\Psi_{sc}^{-2,\infty}}}
\newcommand\Psiscr{\operatorname{\Psi_{sc}^{-2,0}}}
\newcommand\Psiscrm{\operatorname{\Psi_{sc}^{0,2}}}
\newcommand\PsiscHam{\operatorname{\Psi_{sc}^{2,0}}}
\newcommand\Psisci{\operatorname{\Psi_{sc}^{*,*}}}
\newcommand\Psiscid{\operatorname{\Psi_{sc}^{0,0}}}
\newcommand\Psiscis{\operatorname{\Psi_{sc}^{0,\infty}}}
\newcommand\Psiscsi{\operatorname{\Psi_{sc}^{-\infty,0}}}
\newcommand\Psiscs{\operatorname{\Psi_{sc}^{-\infty,\infty}}}
\newcommand\Psiscalg{\operatorname{\Psi_{sc}^{\infty,-\infty}}}
\newcommand\nullHam{{\mathcal N}}
\newcommand\charD{\Sigma_{\Delta-\lambda^2}}
\newcommand\charLap{\Sigma_{\Delta-\lambda}}
\newcommand\Snl{\Sn_{\lambda}}
\newcommand\SNl{\SN_{\lambda}}
\newcommand\gammat{\tilde\gamma}
\newcommand\gammasc{\gamma}
\newcommand\Tau{\mathcal{T}}
\newcommand\taut{\tilde\tau}
\newcommand\taub{\bar\tau}
\newcommand\Nout{N^+_{\lambda}}
\newcommand\Nin{N^-_{\lambda}}
\newcommand\Nio{N^{\pm}_{\lambda}}
\newcommand\El{E_{\lambda}}
\newcommand\Elt{\tilde E_{\lambda}}
\newcommand\Eil{E^i_{\lambda}}
\newcommand\Ejl{E^j_{\lambda}}
\newcommand\Eajl{E^{\alpha_j}_{\lambda}}
\newcommand\Eilt{\tilde E^i_{\lambda}}
\newcommand\Np{N^+}
\newcommand\Nm{N^-}
\newcommand\Npm{N^{\pm}}
\newcommand\Fin{F^-(\lambda)}
\newcommand\Fini{F^-_i(\lambda)}
\newcommand\Fout{F^+(\lambda)}
\newcommand\Fouti{F^+_i(\lambda)}
\newcommand\Foutj{F^+_j(\lambda)}
\newcommand\Rout{R^+_{\lambda}}
\newcommand\Routl{R^+_{\lambda^2}}
\newcommand\Routsgnl{R^{\sign\lambda}_{\lambda^2}}
\newcommand\Rin{R^-_{\lambda}}
\newcommand\Rinl{R^-_{\lambda^2}}
\newcommand\Rinsgnl{R^{-\sign\lambda}_{\lambda^2}}
\newcommand\Rio{R^{\pm}_{\lambda}}
\newcommand\Riol{R^{\pm}_{\lambda^2}}
\newcommand\Roi{R^{\mp}_{\lambda}}
\newcommand\Roil{R^{\mp}_{\lambda^2}}
\newcommand\Riob{R^{\pm}}
\newcommand\Roib{R^{\mp}}
\newcommand\Tio{T^{\pm}}
\newcommand\Tiob{T^{\pm}_{\ff}}
\newcommand\Toi{T^{\mp}}
\newcommand\Toib{T^{\mp}_{\ff}}
\newcommand\TIiob{T_I^{\pm}}
\newcommand\Rinb{R^-}
\newcommand\Rinbsgnl{R^{-\sign\lambda}}
\newcommand\Tin{T^-}
\newcommand\Tinb{T^-_{\ff}}
\newcommand\TIinb{T^-_I}
\newcommand\Routb{R^+}
\newcommand\Routbsgnl{R^{\sign\lambda}}
\newcommand\Tout{T^+}
\newcommand\Toutb{T^+_{\ff}}
\newcommand\TIoutb{T^+_I}
\newcommand\Rlkf{(|\xib|^2-(\lambda-i0)^2)^{-1}}
\newcommand\Rlk{\rho_0(\lambda)}
\newcommand\Rmlk{\rho_0(-\lambda)}
\newcommand\Rpmlk{\rho_0(\pm\lambda)}
\newcommand\Rlka{\rho_1(\lambda)}
\newcommand\Rlkb{\rho_2(\lambda)}
\newcommand\Rilk{\rho_i(\lambda)}
\newcommand\reduced{\natural}
\newcommand\Rlf{R_0(\lambda)}
\newcommand\Rla{R_1(\lambda)}
\newcommand\Rlb{R_2(\lambda)}
\newcommand\Ril{R_i(\lambda)}
\newcommand\Rlj{R_j(\lambda)}
\newcommand\Rlft{R_0(\lambda)}
\newcommand\Rflambda{R_0^{\reduced}(\sigma)}
\newcommand\RV{R^{\reduced}_V}
\newcommand\Rfsigma{R_0^{\reduced}(\sigma)}
\newcommand\Rfsigmah{R_0^{\reduced}(\sigma^{1/2})}
\newcommand\Rfzero{R_0^{\reduced}(0)}
\newcommand\RlV{R^{\reduced}_V(\sigma)}
\newcommand\RlVi{R^{\reduced}_{V_i}(\sigma)}
\newcommand\RlVt{R_V(\lambda)}
\newcommand\RlVtL{{R}_V^L(\lambda)}
\newcommand\RlVtR{{R}_V^R(\lambda)}
\newcommand\RlVit{{R}_{V_i}(\lambda)}
\newcommand\RlVta{{R}_V^{(1)}(\lambda)}
\newcommand\RlVtk{{R}_V^{(k)}(\lambda)}
\newcommand\RlVatV{{R}_{V_{\alpha}}(\lambda)V_{\alpha}}
\newcommand\RlVatVa{{R}_{V_{\alpha_1}}(\lambda)V_{\alpha_1}}
\newcommand\RlVatVb{{R}_{V_{\alpha_2}}(\lambda)V_{\alpha_2}}
\newcommand\RlVatVk{{R}_{V_{\alpha_k}}(\lambda)V_{\alpha_k}}
\newcommand\RlVatVkk{{R}_{V_{\alpha_{k+1}}}(\lambda)V_{\alpha_{k+1}}}
\newcommand\RlVaptV{{R}_{V_{\alpha'}}(\lambda)V_{\alpha'}}
\newcommand\RlVapptV{{R}_{V_{\alpha''}}(\lambda)V_{\alpha''}}
\newcommand\RlVajtV{{R}_{V_{\alpha_j}}(\lambda)V_{\alpha_j}}
\newcommand\RlVaktV{{R}_{V_{\alpha_k}}(\lambda)V_{\alpha_k}}
\newcommand\RlVakktV{{R}_{V_{\alpha_{k+1}}}(\lambda)V_{\alpha_{k+1}}}
\newcommand\Tl{T(\lambda)}
\newcommand\Tlt{\tilde\Tl}
\newcommand\Tltp{\tilde T'(\lambda)}
\newcommand\Tltpp{\tilde T''(\lambda)}
\newcommand\Tli{T_i(\lambda)}
\newcommand\Tlit{\tilde\Tli}
\newcommand\Tlip{T_i'(\lambda)}
\newcommand\Tlipp{T_i''(\lambda)}
\newcommand\Tlj{T_j(\lambda)}
\newcommand\Tla{T_{\alpha}(\lambda)}
\newcommand\Tlaa{T_{\alpha_1}(\lambda)}
\newcommand\Tlab{T_{\alpha_2}(\lambda)}
\newcommand\Tlak{T_{\alpha_k}(\lambda)}
\newcommand\Tlakt{\tilde\Tlak}
\newcommand\Tlaj{T_{\alpha_j}(\lambda)}
\newcommand\Tlajj{T_{\alpha_{j+1}}(\lambda)}
\newcommand\Tlajp{T_{\alpha_j}'(\lambda)}
\newcommand\Tlajpt{\tilde\Tlajp}
\newcommand\Tlajt{\tilde\Tlaj}
\newcommand\Tlakk{T_{\alpha_{k+1}}(\lambda)}
\newcommand\Tlakkp{T_{\alpha_{k+1}}'(\lambda)}
\newcommand\Tlap{T_{\alpha'}(\lambda)}
\newcommand\Tlapt{\tilde\Tlap}
\newcommand\Tlapp{T_{\alpha''}(\lambda)}
\newcommand\Tkl{T^{(k)}(\lambda)}
\newcommand\Tcl{T^{\flat}(\lambda)}
\newcommand\Fl{F(\lambda)}
\newcommand\BlVt{\tilde B_V(\lambda)}
\newcommand\KBlVt{K_{\BlVt}}
\newcommand\BlVaat{B_{V_{\alpha_1}}(\lambda)}
\newcommand\BV{B_V}
\newcommand\Bone{B_1}
\newcommand\Btwo{B_2}
\newcommand\Bthree{B_3}
\newcommand\Banyj{B_j}
\newcommand\PlV{P_V(\lambda)}
\newcommand\PlVc{P_V^{\flat}(\lambda)}
\newcommand\Pl{P_0(\lambda)}
\newcommand\SVl{S_V(\lambda)}
\newcommand\Sjr{S_j^{\reduced}}
\newcommand\Rkp{{\mathcal R}^k_+}
\newcommand\Rkm{{\mathcal R}^k_-}
\newcommand\Rkpm{{\mathcal R}^k_{\pm}}
\newcommand\Phys{{\mathcal P}}
\newcommand\Pc{\overline{\mathcal P}}
\newcommand\pip{\pi^{\perp}}
\newcommand\pipa{\pi_\partial}
\newcommand\gammapa{\gamma_\partial}
\newcommand\pipah{\hat\pi_\partial}
\newcommand\pit{\tilde\pi}
\newcommand\xit{\tilde\xi}
\newcommand\zetat{\tilde\zeta}
\newcommand\etat{\tilde\eta}
\newcommand\sigmat{\tilde\sigma}
\newcommand\sigmahat{\hat\sigma}
\newcommand\thetat{\tilde\theta}
\newcommand\psit{\tilde\psi}
\newcommand\phit{\tilde\phi}
\newcommand\chit{\tilde\chi}
\newcommand\rhot{\tilde\rho}
\newcommand\xib{\bar\xi}
\newcommand\zetab{\bar\zeta}
\newcommand\thetab{\bar\theta}
\newcommand\etab{\bar\eta}
\newcommand\iotal{\iota_{\lambda}}
\newcommand\rhoat{\rhot_{\alpha_1}}
\newcommand\Lambdat{\tilde\Lambda}
\newcommand\Lambdati{\tilde\Lambda^{\text{in}}}
\newcommand\Lambdato{\tilde\Lambda^{\text{out}}}
\newcommand\Lambdatp{\tilde\Lambda^{\text{prop}}}
\newcommand\Lambdai{\Lambda^{\text{in}}}
\newcommand\Lambdao{\Lambda^{\text{out}}}
\newcommand\poles{\Lambda'}
\newcommand\rpoles{\Lambda_p}
\newcommand\thresholds{\Lambda}
\newcommand\Vt{\tilde V}
\newcommand\It{\tilde I}
\newcommand\half{{\frac{1}{2}}}
\newcommand\sigmah{\sigma^{1/2}}
\newcommand\bX{\partial X}
\newcommand\bXb{\partial \Xb}
\newcommand\Deltabt{\tilde\Delta_0}
\newcommand\strip{\Omega_T}
\newcommand\Kf{K^{\flat}}
\newcommand\Gs{G^{\sharp}}
\newcommand\Gt{\tilde G}
\newcommand\Osc{\sci\Omega}
\newcommand\OSc{{}^\Scl\Omega}
\newcommand\Osch{\sci\Omega^{\half}}
\newcommand\Oscmh{\sci\Omega^{-\half}}
\newcommand\Isc{I_{sc}}
\newcommand\os{{\text{os}}}
\newcommand\Qzl{Q^0_{-\lambda}}
\newcommand\Lie{{\mathcal L}}
\newcommand\bl{{\text b}}
\newcommand\scl{{\text{sc}}}
\newcommand\sccl{{\text{scc}}}
\newcommand\Scl{{\text{Sc}}}
\newcommand\ScLl{{\text{Sc,L}}}
\newcommand\ScRl{{\text{Sc,R}}}
\newcommand\Sccl{{\text{Scc}}}
\newcommand\sus{{\text{sus}}}
\newcommand\ssl{{\text{ss}}}
\newcommand\XXb{X^2_\bl}
\newcommand\XXbt{\Xt^2_\bl}
\newcommand\XXsc{X^2_\scl}
\newcommand\XXsct{\Xt^2_\scl}
\newcommand\XXSc{X^2_\Scl}
\newcommand\XXSct{\Xt^2_\Scl}
\newcommand\XXScL{X^2_\ScLl}
\newcommand\XXScR{X^2_\ScRl}
\newcommand\MMsc{M^2_\scl}
\newcommand\Deltab{\Delta_\bl}
\newcommand\Deltasc{\Delta_\scl}
\newcommand\DeltaSc{\Delta_\Scl}
\newcommand\DeltaScL{\Delta_\ScLl}
\newcommand\DeltaScR{\Delta_\ScRl}
\newcommand\prs{\sigma}
\newcommand\Nsc{N_\scl}
\newcommand\Nscp{N_{\scl,p}}
\newcommand\Nff{N_{\ff}}
\newcommand\Nffz{N_{\ff,0}}
\newcommand\Nffzp{N_{\ff,0,p}}
\newcommand\Nffl{N_{\ff,l}}
\newcommand\Nffml{N_{\ff,-l}}
\newcommand\Nmf{N_{\mf}}
\newcommand\Nmfz{N_{\mf,0}}
\newcommand\Nmfl{N_{\mf,l}}
\newcommand\Nmfml{N_{\mf,-l}}
\newcommand\ffb{\operatorname{bf}}
\newcommand\Ffb{\operatorname{bf'}}
\newcommand\ffsc{\operatorname{sf}}
\newcommand\ffSc{\operatorname{sf_C}}
\newcommand\Ffsc{\operatorname{sf'}}
\newcommand\rff{\rho_{\ff}}
\newcommand\rmf{\rho_{\mf}}
\newcommand\rffb{\rho_{\ffb}}
\newcommand\rffsc{\rho_{\ffsc}}
\newcommand\rFfsc{\rho_{\Ffsc}}
\newcommand\rffSc{\rho_{\ffSc}}
\newcommand\rinf{\rho_{\infty}}
\newcommand\CL{C_L}
\newcommand\CR{C_R}
\newcommand\betab{\beta_\bl}
\newcommand\betasc{\beta_\scl}
\newcommand\betaSc{\beta_\Scl}
\newcommand\BetaSc{\bar\beta_\Scl}
\newcommand\betaScL{\beta_\ScLl}
\newcommand\betaScR{\beta_\ScRl}
\newcommand\ScT{{}^\Scl T^*}
\newcommand\SccT{{}^\Scl \bar T^*}
\newcommand\ScS{{}^\Scl S^*}
\newcommand\Tb{{}^\bl T}
\newcommand\Tsc{{}^\scl T}
\newcommand\TSc{{}^\Scl T}
\newcommand\CSc{C_\Scl}
\newcommand\Lambdasc{{}^\scl\Lambda}
\newcommand\XXXb{X^3_\bl}
\newcommand\XXXsc{X^3_\scl}
\newcommand\XXXSc{X^3_\Scl}
\newcommand\XXXScO{X^3_{\Scl,O}}
\newcommand\XXXScF{X^3_{\Scl,F}}
\newcommand\XXXScS{X^3_{\Scl,S}}
\newcommand\XXXScC{X^3_{\Scl,C}}
\newcommand\KDsc{\operatorname{KD^{\half}_\scl}}
\newcommand\KDSc{\operatorname{KD^{\half}_\Scl}}
\newcommand\KDScEF{\operatorname{KD^{E,F}_\Scl}}
\newcommand\Oh{\operatorname{\Omega^{\half}}}
\newcommand\WFSc{\WF_\Scl}
\newcommand\WFtSc{\WF_{\text 3sc}}
\newcommand\WFScmf{\WF_{\Scl,\mf}}
\newcommand\WFScff{\WF_{\Scl,\ff}}
\newcommand\WFScs{\WF_{\Scl,\prs}}
\newcommand\WFScp{\WF'_\Scl}
\newcommand\WFScmfp{\WF'_{\Scl,\mf}}
\newcommand\WFScffp{\WF'_{\Scl,\ff}}
\newcommand\WFScsp{\WF'_{\Scl,\prs}}
\newcommand\Diffscc{\Diff_\sccl}
\newcommand\DiffSc{\Diff_\Scl}
\newcommand\DiffScc{\Diff_\Sccl}
\newcommand\DiffscI{\Diff_{\scl,\text{I}}}
\newcommand\VscI{\Vf_{\scl,\text{I}}}
\newcommand\DiffsV{\operatorname{Diff}_{\sus(V)}}
\newcommand\DiffsVsc{\operatorname{Diff}_{\sus(V),\scl}}
\newcommand\DiffsVCsc{\operatorname{Diff}_{\sus(V)-C,\scl}}   
\newcommand\Psisc{\Psop_\scl}
\newcommand\Psiscc{\Psop_\sccl}
\newcommand\Psiss{\Psop_\ssl}
\newcommand\Psisch{\Psop_{\scl,h}}
\newcommand\Psiscch{\Psop_{\sccl,h}}
\newcommand\PsiSc{\Psop_\Scl}
\newcommand\PsiScph{\Psop_{\Scl,\phi}}
\newcommand\PsiScra{\Psop_{\Scl,\rho^\sharp_a}}
\newcommand\PsiScc{\Psop_\Sccl}
\newcommand\PsiSccml{\Psop^{m,l}_\Sccl}
\newcommand\PsiScxx{\Psop^{*,*}_\Scl}
\newcommand\PsiScml{\Psop^{m,l}_\Scl}
\newcommand\PsiScmz{\Psop^{m,0}_\Scl}
\newcommand\PsiScmmz{\Psop^{-m,0}_\Scl}
\newcommand\PsiSckz{\Psop^{k,0}_\Scl}
\newcommand\PsiScmmml{\Psop^{-m,-l}_\Scl}
\newcommand\Psiscmkk{\Psop^{-k,k}_\scl}
\newcommand\Psiscmmmkk{\Psop^{-m-k,k}_\scl}
\newcommand\Psiscmoo{\Psop^{-1,1}_\scl}
\newcommand\Psiscmz{\Psop^{m,0}_\scl}
\newcommand\Psiscmmz{\Psop^{-m,0}_\scl}
\newcommand\PsiSckmkl{\Psop^{km,kl}_\Scl}
\newcommand\PsiScmplp{\Psop^{m',l'}_\Scl}
\newcommand\PsiScmmpllp{\Psop^{m+m',l+l'}_\Scl}
\newcommand\Psiscml{\Psop^{m,l}_\scl}
\newcommand\PsiScid{\Psop^{0,0}_\Scl}
\newcommand\PsiSczo{\Psop^{0,1}_\Scl}
\newcommand\PsiScmii{\Psop^{-\infty,\infty}_\Scl}
\newcommand\PsiScmiz{\Psop^{-\infty,0}_\Scl}
\newcommand\PsiScmoo{\Psop^{-1,1}_\Scl}
\newcommand\PsisCid{\Psop^{0,0}_{\scl-C}}
\newcommand\PsisC{\Psop_{\scl-C}}
\newcommand\Psiinf{\Psop_{\infty}}
\newcommand\Psiinfid{\Psop_{\infty}^0}
\newcommand\PsiFinf{\Psop_{\infty-\Fr}}
\newcommand\PsisVscml{\Psop^{m,l}_{\sus(V),\scl}}
\newcommand\PsisVsc{\Psop_{\sus(V),\scl}}
\newcommand\PsisVpsc{\Psop_{\sus(V_p),\scl}}
\newcommand\PsisVCSc{\Psop_{\sus(V)-C,\scl}}
\newcommand\SFinf{S_{\infty-\Fr}}
\newcommand\YsVC{Y^2_{\sus(V)-C,\scl}}
\newcommand\ffYsc{\ffsc_{\sus(V)}}
\newcommand\SXC{S(X;C)}
\newcommand\Ios{I_{\text{os}}}
\newcommand\pbL{\pi^2_{\bl,{\text L}}}
\newcommand\pbR{\pi^2_{\bl,{\text R}}}
\newcommand\pscL{\pi^2_{\scl,{\text L}}}
\newcommand\pscR{\pi^2_{\scl,{\text R}}}
\newcommand\PbO{\pi^3_{\bl,{\text O}}}
\newcommand\PscO{\pi^3_{\scl,{\text O}}}
\newcommand\PScO{\pi^3_{\Scl,{\text O}}}
\newcommand\PScF{\pi^3_{\Scl,{\text F}}}
\newcommand\PScC{\pi^3_{\Scl,{\text C}}}
\newcommand\PScS{\pi^3_{\Scl,{\text S}}}
\newcommand\pScL{\pi^2_{\Scl,{\text L}}}
\newcommand\pScR{\pi^2_{\Scl,{\text R}}}
\newcommand\CLF{\CL^F}
\newcommand\CLO{\CL^O}
\newcommand\CLS{\CL^S}
\newcommand\CLC{\CL^C}
\newcommand\DeltaYb{\Delta_{\bl,Y}}
\newcommand\DeltaYsc{\Delta_{\sus-\scl}}
\newcommand\diag{\operatorname{diag}}
\newcommand\Vf{{\mathcal V}}
\newcommand\Vb{{\mathcal V}_{\bl}}
\newcommand\Vsc{{\mathcal V}_{\scl}}
\newcommand\VSc{{\mathcal V}_{\Scl}}
\newcommand\VfI{\Vf_{\text{I}}}
\newcommand\VfIq{\Vf_{\text{I},q}}
\newcommand\scH{{}^\scl H}
\newcommand\scHg{\scH_g}
\newcommand\Hss{H_\ssl}
\newcommand\xh{\hat x}
\newcommand\Yh{\hat Y}
\newcommand\Zh{\hat Z}
\newcommand\Yb{\bar Y}
\newcommand\hb{\bar h}
\newcommand\xih{\hat\xi}
\newcommand\etah{\hat\eta}
\newcommand\muh{\hat\mu}
\newcommand\mub{\bar\mu}
\newcommand\nub{\bar\nu}
\newcommand\mubh{\widehat{\bar\mu}}
\newcommand\yb{\bar y}
\newcommand\zb{\bar z}
\newcommand\ub{\bar u}
\newcommand\Qb{\bar Q}
\newcommand\Wbp{{\bar W}^\perp}
\newcommand\Wp{W^\perp}
\newcommand\Kt{\tilde K}
\newcommand\Wt{\tilde W}
\newcommand\Ut{\tilde U}
\newcommand\yt{\tilde y}
\newcommand\ut{\tilde u}
\newcommand\vt{\tilde v}
\newcommand\ft{\tilde f}
\newcommand\htil{\tilde h}
\newcommand\St{\tilde S}
\newcommand\Pt{\tilde P}
\newcommand\Rt{\tilde R}
\newcommand\qt{\tilde q}
\newcommand\Qt{\tilde Q}
\newcommand\Xb{\bar X}
\newcommand\lambdat{\tilde\lambda}
\newcommand\betat{\tilde\beta}
\newcommand\epst{\tilde\epsilon}
\newcommand\ep{\epsilon}
\newcommand\bt{\tilde b}
\newcommand\Xt{\tilde X}
\newcommand\Mt{\tilde M}
\newcommand\At{\tilde A}
\newcommand\Et{\tilde E}
\newcommand\Ht{\tilde H}
\newcommand\at{\tilde a}
\newcommand\Ct{\tilde C}
\newcommand\pih{\hat\pi}
\newcommand\Rh{\hat R}
\newcommand\Ah{\hat A}
\newcommand\Bh{\hat B}
\newcommand\Ch{\hat C}
\newcommand\Gh{\hat G}
\newcommand\Hh{\hat H}
\newcommand\Qh{\hat Q}
\newcommand\Ph{\hat P}
\newcommand\Nh{\hat N}
\newcommand\Sh{\hat S}
\newcommand\Gcal{{\mathcal G}}
\newcommand\GcalC{{\mathcal G}_C}
\newcommand\Jcal{{\mathcal J}}
\newcommand\JcalC{{\mathcal J}_C}
\setcounter{secnumdepth}{3}
\newtheorem{lemma}{Lemma}[section]
\newtheorem{prop}[lemma]{Proposition}
\newtheorem{thm}[lemma]{Theorem}
\newtheorem{cor}[lemma]{Corollary}
\newtheorem{result}[lemma]{Result}
\newtheorem*{thm*}{Theorem}
\newtheorem*{prop*}{Proposition}
\newtheorem*{conj*}{Conjecture}
\numberwithin{equation}{section}
\theoremstyle{remark}
\newtheorem{rem}[lemma]{Remark}
\theoremstyle{definition}
\newtheorem{Def}[lemma]{Definition}
\newtheorem*{Def*}{Definition}
\def\signature#1#2{\par\noindent#1\dotfill\null\\*
{\raggedleft #2\par}}

\renewcommand{\theenumi}{\roman{enumi}}
\renewcommand{\labelenumi}{(\theenumi)}
\newcommand{\CC}{\mathbb C}
\newcommand{\HHk}{{\mathbb H}^{k+1}}
\newcommand{\HH}{\mathbb H}
\newcommand{\RR}{\mathbb R}
\newcommand{\del}{\partial}
\newcommand{\e}{\epsilon}
\newcommand{\cC}{{\mathcal C}}
\newcommand{\Lap}{\Delta}
\newcommand{\Ree}{\mbox{\rm Re}\,}
\newcommand{\Imm}{\mbox{\rm Im}\,}
\newcommand{\bM}{\overline{M}}
\renewcommand{\bX}{\overline{X}}
\newcommand{\bZ}{\overline{Z}}
\newcommand{\res}{{\mathrm{res}}}
\newcommand{\well}{w^{(\ell)}}

\title{Resolvents and Martin boundaries of product spaces}
\author[Rafe Mazzeo and Andras Vasy]{Rafe Mazzeo and Andr\'as Vasy}
\date{November 29, 2000}
\thanks{R.\ M.\ partially supported by NSF grant \#DMS-99-1975;
A.\ V.\ partially supported by NSF grant \#DMS-99-70607.}
\address{R.\ M.: Department of Mathematics, Stanford University, Stanford,
CA 94305}
\email{mazzeo@math.stanford.edu}
\address{A.\ V.: Department of Mathematics, Massachusetts Institute of
Technology, MA 02139}
\email{andras@math.mit.edu}

\maketitle

\section{Introduction}
Geometric scattering theory, as espoused in \cite{RBMGeo}, is the
study of natural operators such as the Laplacian on non-compact
Riemannian manifolds $X$ with controlled asymptotic geometry using
geometrically-informed, fully (i.e.\ to the extent it is possible)
microlocal methods. The goals of this
subject include the construction of an analytically useful compactification
of $X$ and the definition of an appropriate class of pseudodifferential
operators containing the resolvent of the Laplacian, the Schwartz kernel of
which is a particularly simple distribution on this compactification. In
this paper we examine the case where $X$ is a product of asymptotically
hyperbolic (or conformally compact, as they are often called) spaces from
this point of view. This is intended both as an initial application of these
methods to higher rank symmetric spaces and their geometric generalizations,
and also as a relatively simple (although still surprisingly complicated)
example which should provide a guide for what to expect in further
development in this area.

The results here tend to be notationally complicated, so for the purposes
of this introduction we state a simple, yet representative, result
concerning the asymptotics of the resolvent applied to a Schwartz function.
Suppose that $(M_j,g_j)$, $j=1,2$, are (conformally compact) asymptotically
hyperbolic, $\dim M_j= k_j+1$, $H_j=-\Delta_{g_j}$ positive,
$\phi_{ji}$ are the $L^2$
eigenfunctions of $H_j$ with eigenvalue $\lambda_{ji} <  k_j^2/4$. The usual
compactification of $M_j$ as a $\Cinf$ manifold with boundary is denoted
$\bM_j$, and has boundary defining function $x_j$. Adjoining $\rho_j=
-1/\log x_j$ to the smooth structure of $\bM_j$ yields the logarithmic blow up
$(\bM_j)_{\log}$.  Let $H$ be the sum of the Laplacians from the two
factors: $H=H_1\otimes\Id+\Id\otimes H_2$, and $R(\mu)$ the resolvent,
$(H - \mu)^{-1}$. Define
\begin{equation}\label{eq:def-Xt}
\Xt=[(\bM_1)_{\log}\times(\bM_2)_{\log};\pa(\bM_1)_{\log}
\times\pa(\bM_2)_{\log}],
\end{equation}
which is a resolution of $\bX=\bM_1\times\bM_2$. By definition,
each $\rho_j=-1/\log x_j$ is smooth on $(\bM_j)_{\log}$ and thus
$\rho$, where $\rho^{-1}=\sqrt{\rho_1^{-2}+\rho_2^{-2}}$ is smooth on $\Xt$.

\begin{thm*} (See Theorem~\ref{thm:cont-spec-asymp}.)
Suppose $f\in\dCinf(\bX)$, $\Real \ni \mu>k^2/4=(k_1^2+k_2^2)/4$.
Then on $\Xt$:
\begin{equation}\begin{split}
& R(\mu-i0)f=x_1^{k_1/2}x_2^{k_2/2}\exp (-i\sqrt{\mu-k^2/4}/\rho)g\\
+ & \sum_{i=1}^{N_1}x_2^{k_2/2+
i\sqrt{\mu-\lambda_{1i}-k_2^2/4}}(\phi_{1i} \otimes g_{1i})
+\sum_{i=1}^{N_2}x_1^{k_1/2+
i\sqrt{\mu-\lambda_{1i}-k_2^2/4}}(g_{2i}\otimes \phi_{2i}),
\end{split}\end{equation}
$g$ polyhomogeneous on $\Xt$, $g_{1i}$ and $g_{2i}$ are polyhomogeneous
on $\bar M_1$ and $\bar M_2$ respectively.
\end{thm*}

The leading term of each part in the asymptotics can be described explicitly.
Since $\phi_{ji}$ decays like $x_j^{k_j/2+\sqrt{k_1^2/4-\lambda_{1i}}}$,
the first term dominates the
other two in the interior of the front face of the blow-up
\eqref{eq:def-Xt}; the second term, involving the eigenfunctions
of $H_1$, is only comparable to it at the lift of $\bM_1\times\pa\bM_2$, while
the third term, involving the eigenfunctions
of $H_2$, is only comparable to it at the lift of $\pa\bM_1\times\bM_2$.

Similar results are valid outside the continuous spectrum, where the
exponents are not pure imaginary, but care needs
to be taken as to which eigenvalue terms appear in the asymptotics.
Moreover, additional singularities appear along the front face
of $\Xt$ when $\mu$ is real, below the spectrum of $H$, which, however,
can be resolved by further real blow-ups (or understood as a type of Legendre
singularity).
We refer to Theorem~\ref{thm:res-set-asymp} for the detailed statement of
the results in that case.
We deduce similar results for the resolvent kernel, which
we use in turn to analyze the Martin compactification
${\overline X}_M$of $M_1\times M_2$.
While the latter behaves (nearly) as expected when $H_1$ and $H_2$ have
no $L^2$ eigenvalues, it experiences a substantial collapse in the presence
of such eigenvalues.

\begin{thm*}
There is a natural continuous surjection $\Xt\to{\overline X}_M$. Moreover,
the following hold.

\begin{enumerate}
\item
If neither $H_1$ nor $H_2$ have $L^2$ eigenvalues,
the restriction of this map to a
neighborhood of the front face of the blow-up in \eqref{eq:def-Xt} is
injective. In general, its injectivity on the `side faces' of $\Xt$
depends on properties of the spherical functions, i.e.\ on
$\frac{d}{d\tau}|_{\tau=0}
R_j(k_j^2/4+\tau^2)$.

\item
If both $H_1$ and $H_2$ have $L^2$ eigenvalues, the Martin boundary
$\pa{\overline X}_M$
is of the form $\pa \bM_1\cup\pa\bM_2\cup\pa \bM_1\times\pa\bM_2\times I$,
$I$ an open interval. This set is naturally identified with a collapsed
version of $\Xt$, and the map $\pa\Xt\to\pa{\overline X}_M$ factors
through it.
\end{enumerate}
\end{thm*}

There is an extensive literature on scattering on conformally compact
manifolds, including especially the important special case of convex cocompact
(and geometrically finite) hyperbolic manifolds. This contains geometric
constructions of the resolvent, scattering matrix and generalized
eigenfunctions, as well as trace formul\ae\ and asymptotics for the counting
function for resonances. For brevity, we list only
\cite{Mazzeo-Melrose:Meromorphic}, \cite{Mazzeo:Hodge}, \cite{Mazzeo:Unique},
\cite{Perry:Laplace-I}, \cite{Hislop:Geometry}, \cite{Zworski:Dimension} for
references. A parallel theory for complex hyperbolic manifolds and their
perturbations is initiated in \cite{Epstein-Melrose-Mendoza:Resolvent},
and while not appearing explicitly, this theory extends to manifolds with the
asymptotic structure of quaternionic hyperbolic spaces and the Cayley plane.

Concerning higher rank noncompact symmetric spaces (and in less generality,
their quotients too), the compactification theory is now well-understood
\cite{Guivarch-Ji-Taylor:Compactifications}, as well as at least some
aspects of the analysis of the Laplacian. However, the current methods here
rely heavily on the special algebraic structure of these spaces, so
extensions of these results to even relatively modest geometric perturbations
of these spaces are basically not understood at all.

Alongside this is the study of the asymptotically Euclidean scattering
metrics initiated in \cite{RBMSpec}, \cite{RBMZw}. This theory extends the
extensive classical literature, but has led to a new and detailed
understanding of quantum N-body scattering where the geometry has much in
common with that of flats in non-compact symmetric spaces, see
e.g.\ \cite{Vasy:Structure}, \cite{Hassell:Plane}, \cite{Vasy:Propagation-2}
and \cite{Vasy:Bound-States}.

These different geometric settings have led to the understanding of diverse
analytic phenomena which may sometimes be traced to the effects of the
asymptotically flat or negative curvature.  One of the attractions of
studying higher rank symmetric spaces from the point of view of geometric
scattering theory is to isolate the specific ways in which the flat and
negatively curved directions interact with one another and affect the
analysis.  The simplest setting where these sorts of effects might be seen
is on the product $\bX = \bM_1 \times \bM_2$, where both factors $(\bM_j,g_j)$
are conformally compact metrics. Recall that $(\bM,g)$ is conformally compact
if $\bM$ is a smooth compact manifold with boundary, such that for some
defining function $x$ for $\del \bM$, $g$ takes the form $(dx^2 + h)/x^2$,
where $h$ is a nonnegative smooth symmetric $2$-tensor which restricts to a
nondegenerate metric on $\del \bM$. (Note, however, that only the conformal
class of $h$ on $\del \bM$ is well-defined from $g$.) The prototype of a
conformally compact manifold is hyperbolic space, and so the prototypes for
the spaces we consider here are products of hyperbolic spaces. The reader
should be aware, though, that $\bM_1 \times \bM_2$ is asymptotically like
the product of hyperbolic spaces $\HH^{n_1} \times \HH^{n_2}$ ($n_i =
\dim M_i$) only near $\del \bM_1 \times \del \bM_2$, but elsewhere may be
a rather severe metric and topological perturbation of this model space.

The basic questions we ask here concerning the product space
$X$ with metric $g = g_1 + g_2$ concern the resolvent $R_X(\lambda) =
(\Lap_g - \lambda)^{-1}$ of its Laplacian $\Lap_g$. We sometimes
write $R_X(\lambda)$ as $R_g(\lambda)$ or simply $R(\lambda)$. This is a
holomorphic family of elliptic pseudodifferential operators of order $-2$
in the resolvent set $\CC \setminus \spec(\Lap_g)$, which contains at least
the complement of the positive real axis. The precise extent of the spectrum
is straightforward to determine from the spectra on either factor, but the
first substantial question is to understand the behaviour of $R(\lambda)$
as $\lambda$ approaches the (continuous part of the) spectrum. Existence of
a limit, in an appropriate sense, is known as the limiting absorption
principle. Even better is the existence of a meromorphic continuation of
$R(\lambda)$ beyond the spectrum. Usually, when it exists, this continuation
lives on some Riemann surface covering the complex plane and ramified at
the thresholds of the spectrum. Meromorphic continuations of this type are
known to exist for the Laplacian of conformally compact manifolds
\cite{Mazzeo-Melrose:Meromorphic}. Further important questions concern the
geometric space which is a resolution of $X \times X$ obtained by a process
of real blow-up, and on which the Schwartz kernel of the resolvent lives as a
polyhomogeneous (or more generally, a Legendrian) distribution. Detailed
knowledge of this space is central in understanding the finer properties of
the resolvent, and conversely is key in its initial geometric construction.
Geometric and analytic compactifications of the space $X$ itself are quite
relevant to this. There are two well-known compactifications which are
closely related to our methods: the geodesic (also known as the conic)
compactification, and the Martin compactification, which is defined using
function theory, specifically the space of positive solutions of
$(\Lap_g - \lambda)u = 0$ for $\lambda$ real and below the bottom of
$\spec(\Lap_g)$. This Martin compactification is known in many instances,
including for the product of hyperbolic spaces \cite{Giulini-Woess:Martin}
and more recently for general symmetric spaces of noncompact type
\cite{Guivarch-Ji-Taylor:Compactifications}. The `resolvent compactification'
of $X \times X$ we construct below, and on which the Schwartz kernel of the
resolvent lives, has a structure on some of its hypersurface boundary
faces which is inherited from the geodesic and Martin compactifications.

These questions together have led to the somewhat modest goals of the present
paper. In the next section we present a contour integral formula for
$R_X(\lambda)$ in terms of the resolvents of the two factors $R_{M_j}
(\lambda)$. A related expression, written as an integral over the spectral
measure, derived in the context in Euclidean scattering, appears in work of
Ben-Artzi and Devinatz \cite{Ben-Artzi-Devinatz:Resolvent}. The representation
formula here holds in great generality. In \S 3 we specialize and review
the detailed structure of the resolvent $R_{M}$ when $M$ is a conformally
compact manifold; some auxiliary estimates required later are also
derived here. An immediate consequence of the representation formula of \S 2
is the existence of an appropriate meromorphic continuation for the
resolvent $R_X(\lambda)$, and we describe this in \S 4. After this,
\S 5 contains a discussion of general compactification theory, specialized
to this context. The main work is done in \S6 and \S 7, where we describe
the asymptotics of $R_X(\mu)f$ when $f$ is Schwartz, first when $\mu$ is
in the resolvent set and then when $\mu$ is in the main sheet of continuous
spectrum of $H$. The first application of this is given in \S 8, where
we construct the `resolvent double space', a resolution of $X \times X$
which carries the Schwartz kernel of $R_X(\mu)$ in as simple a fashion
as possible. Finally, in \S 9, we use the resolvent asymptotics to
determine the Martin compactification of $X$.

There are several features of this work to which we wish to draw particular
attention. First, the main tool in deriving the resolvent asymptotics
is stationary phase, which is in many respects local. The previous
identification of the Martin compactification (at least for the product
of hyperbolic spaces) in \cite{Giulini-Woess:Martin} relies heavily
on global heat kernel bounds, which we feel are intrinsically more
complicated. Note also that by using stationary phase we are taking
advantage of the oscillatory nature of the resolvent, even when studying
it for certain values of $\mu$ where it has been more traditional to
rely on `positivity methods' such as the maximum principle, the Harnack
inequality, etc. Another interesting feature here is the surprisingly
complicated way the existence of bound states, i.e. $L^2$ eigenvalues,
for the Laplacian on either factor affects the asymptotics and the
structure of the Martin boundary. Finally, we have given a detailed
description of the smooth structure of the various compactifications
we construct; this aspect is usually neglected in other discussions of
compactification theory, but as we show, plays a significant role.

Our intention is that the investigations here will form the basis
for a more thorough investigation of the geometric scattering theory
for higher rank spaces.

The authors wish to thank Richard Melrose for helpful advice, and
also Lizhen Ji for encouraging us to study the Martin compactification.

\section{Resolvent formula}
Let $H_1$, $H_2$ be self-adjoint operators on Hilbert spaces $V_1$ and
$V_2$ which are bounded below; the precise structure of their spectra
will be unimportant for the present. We denote $\inf \spec (H_j) =
\lambda_0(H_j)$. Now let
$$
H=H_1\otimes\Id+\Id\otimes H_2
$$
be the self-adjoint operator on the completed tensor product space
$V_1 \hat \otimes V_2$. Then $\inf\spec(H) \equiv \lambda_0(H) =
\lambda_0(H_1) + \lambda_0(H_2)$. Let
\[
R_j(\mu_j) = (H_j - \mu_j)^{-1}, \quad j = 1,2, \qquad \mbox{and}\quad
R(\mu) = (H - \mu)^{-1}
\]
be the resolvents of the $H_j$ and $H$, respectively. This setting has been
investigated by Ben-Artzi and Devinatz in \cite{Ben-Artzi-Devinatz:Resolvent},
mostly from the point of view of the limiting absorption princple, i.e.\ the
existence of the boundary values $R(\mu\pm i0)$, $\mu\in[\lambda_0(H_1)+
\lambda_0(H_2),\infty)$ on suitable weighted spaces, under the assumption
that the limiting absorption principle holds for the $R_j$ individually.

One of our first goals is to show that if both $R_j$ admit meromorphic
continuations, then $R(\mu)$ does as well. Later we wish to obtain precise
asymptotics for the Schwartz kernel of $R(\mu)$ for $\mu$ in the
resolvent set and in the continuous spectrum. In this section we derive a
representation of $R(\mu)$ as a
contour integral which will be useful for both of these purposes.

\begin{figure}[ht]
\begin{center}
\mbox{\epsfig{file=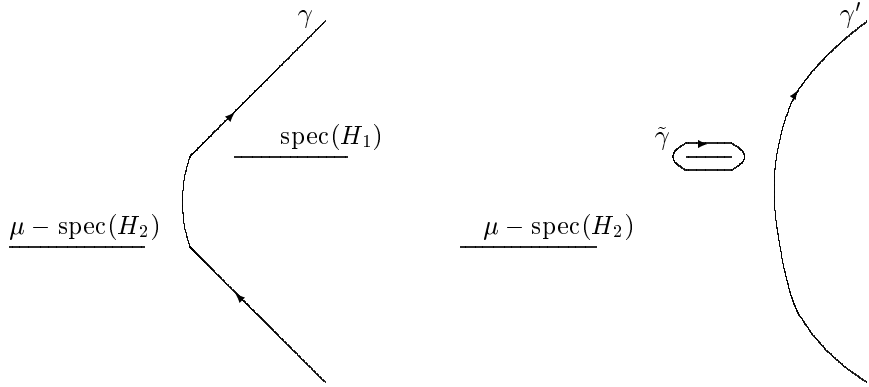}}
\end{center}
\caption{Contours of integration in the $\mu_1$ plane.}
\label{fig:contour}
\end{figure}

So, fix $\mu$ in the resolvent set $\CC \setminus [\lambda_0(H_1)+ \lambda_0
(H_2),\infty)$, and let $\gamma$ be a parametrized curve in the complex plane
which is disjoint from $[\lambda_0(H_1),\infty)$ and such that $\mu-\gamma$
is disjoint from $[\lambda_0(H_2),\infty)$, or in other words, such that
$\gamma$ does not intersect $[\lambda_0(H_1),\infty)\cup(\mu-[\lambda_0(H_2),
\infty))$. Suppose also that $\gamma(t)=c_\pm t$ for $\pm t\geq T>0$ with
$\im c_\pm > 0$. This is illustrated in Figure~\ref{fig:contour}. The precise
values of $c_\pm$ are unimportant for our purposes, and indeed there is
considerably more leeway than this in choosing $\gamma$, but the definite
separation of $\gamma$ from the spectra, ensured by $\im c_\pm>0$, is
crucial. (One can allow slighly subconic separation, but this is of no
interest here.) Then we claim that
\begin{equation}
R(\mu)=\frac{1}{2\pi i}\int_{\gamma} R_1(\mu_1)
R_2(\mu-\mu_1)\,d\mu_1.
\label{eq:resform}
\end{equation}

To prove this formula, first observe that since
$$
\|R_j(\mu_j)\|\leq |\im\mu_j|^{-1},
$$
the norm of the integrand in (\ref{eq:resform}) (as a bounded operator on
$L^2$) is estimated by $C(1+|t|)^{-2}$, and hence the integral converges.
Next, note that it suffices to show that both sides of (\ref{eq:resform})
produce the same result when restricted to the range of $\chi_I(H_1)\otimes
\Id$ where $I$ is any compact interval and $\chi_I$ its characteristic
function, because the union of these ranges is dense.

Fixing the interval $I$, the integrand $R_1(\mu_1)R_2(\mu - \mu_2)$
is holomorphic for $\mu_1 \notin I\cup(\mu-\spec(H_2))$. Therefore we may
deform the contour to one which is the union of two curves, $\hat{\gamma}$
and $\gammat$, where $\hat{\gamma}$ agrees with $\gamma$ for $t\geq T'>0$
and intersects the real axis precisely once, somewhere in $(\sup I,\infty)$,
while $\gammat$ surrounds $I$ once. Thus
\begin{equation}\begin{split}
\frac{1}{2\pi i}&\int_{\gamma} R_1(\mu_1)
R_2(\mu-\mu_1)\,d\mu_1\\
&= \frac{1}{2\pi i}\int_{\hat{\gamma}} R_1(\mu_1)
R_2(\mu-\mu_1)\,d\mu_1+ \frac{1}{2\pi i}\int_{\gammat} R_1(\mu_1)
R_2(\mu-\mu_1)\,d\mu_1.
\end{split}\end{equation}

But now, if $\hat{\gamma}$ is moved to infinity, the first integral on the
right tends to 0, as follows by directly estimating the integral. On the
other hand, letting $\gammat$ tend to $I$ and applying Stone's theorem, the
second integral on the right is the same as
\[
\int_I R_2(\mu-\mu_1)\,dE_1(\mu_1).
\]
This is identical to $R(\mu)$ on the range of $\chi_I(H_1) \otimes \Id$
since applying $H-\mu$ to it gives (with a slight abuse of notation)
\[
\int_I (H_1 - \mu_1 + H_2 - (\mu - \mu_1)) R_2(\mu - \mu_1)\, dE_1(\mu_1).
\]
\[
= \int_I ((H_1-\mu_1)R_2(\mu-\mu_1)+\Id)\,dE_1(\mu_1)
=\int_I\,dE_1(\mu_1) = \chi_I(H_1)\otimes\Id,
\]
i.e.\ the identity on this subspace. Thus (\ref{eq:resform}) is established.

\section{Resolvents of conformally compact manifolds}\label{sec:conf-compact}
In this section we briefly collect some facts about the resolvent family
$R_M(\lambda)$, which we abbreviate simply as $R(\lambda)$ for the
duration of this section, when $(M,g)$ is a conformally compact manifold.
These are all discussed and proved in \cite{Mazzeo-Melrose:Meromorphic},
to which we refer for all details.

Recall that $M$ is identified with the interior of a compact smooth manifold
with boundary $\bM$. The compactification $\bM$ is geometrically natural:
$\bM$ may also be identified with the both the geodesic and Martin
compactifications, as we discuss further below. Locally, $R(\lambda)$ is a
pseudodifferential operator of order $-2$, but our focus is on understanding
the behavior of its Schwartz kernel $K(z,z')= K(\lambda;z,z')$, $z, z' \in M$,
when one or both of the variables tend to infinity in $M$, i.e. to a point of
$\del M$. This can occur in various ways, of course, and the most efficient
way to encode this information is to consider $K$ as a distribution on the
$0$-stretched product $\bM_0^2$ introduced in
\cite{Mazzeo-Melrose:Meromorphic} and \cite{Mazzeo:Hodge}. This space is
obtained from $\bM^2$ by blowing up the boundary of the diagonal $\del \Delta
\iota = \del \{z = z'\}$; equivalently, $\bM^2_0 = [\bM^2; \del \Delta
\iota]$ is the disjoint union of $\bM^2 \setminus \del \Delta \iota $ with
the interior spherical normal bundle of $\del \Delta \iota$ at the corner
$\del \bM \times \del \bM$, this set then
being endowed with the minimal ${\mathcal C}^\infty$ structure containing
the lifts of smooth functions on $\bM^2$ and polar coordinates around
$\del \Delta \iota$.

$\bM^2_0$ has three different boundary hypersurface faces: the left
face which covers (indeed, is identified with) $\del \bM \times \bM$, the
right face which covers $\bM \times \del \bM$, and the new front face
covering $\del \Delta \iota$. We denote these $B_{10}$, $B_{01}$ and $B_{11}$,
and their boundary defining functions $\rho_{10}$, $\rho_{01}$ and
$\rho_{11}$, respectively. Writing $z = (x,y)$ and $z' = (x',y')$ near
$\del \bM$, then
\[
\rho_{11} = \sqrt{x^2 + (x')^2 + |y-y'|^2},
\qquad \mbox{and} \qquad
\rho_{10} = x/\rho_{11}, \quad \rho_{01} = x'/\rho_{11}.
\]

If $\dim M = n = k+1$, then write $\lambda = \zeta (k - \zeta)$, where by
convention the region $\Ree \zeta > k/2$ corresponds to the resolvent set
$\CC \setminus [k^2/4,\infty)$ of $\Lap_g$. Thus for $\lambda$ in the
resolvent set,
$$
\zeta=k/2+i\sqrt{\lambda-k^2/4},
$$
where the branch of the square root is chosen so that its imaginary part
is negative.

\begin{thm} (Mazzeo and Melrose,\
\cite[Theorem~7.1]{Mazzeo-Melrose:Meromorphic})
For $\zeta$ in the half-plane $\Ree \zeta > k/2$, the Schwartz kernel
$K_\zeta$ of $R(\zeta(k - \zeta))$ is a polyhomogeneous distribution on
$\bM^2_0$ with the following properties. First, $R(\zeta(k - \zeta))$ has a decomposition
$$
R(\zeta(k - \zeta))=R'(\zeta(k - \zeta))+R''(\zeta(k - \zeta)),
$$
where $R'$ is an element in the small calculus $\Psi_0^{-2}(\bM)$
of $0$-pseudodifferential operators on $M$ and $R''$ is a residual
element in the large calculus $\Psi_0^{-\infty,\zeta,\zeta}(\bM)$ of
$0$-pseudodifferential operators. This means that the Schwartz kernel of
$R'$ has a standard polyhomogeneous singularity corresponding to
pseudodifferential order $-2$ at the lifted diagonal $\Delta \iota_0$
in $\bM_0^2$ and vanishes to infinite order along $B_{10}$ and $B_{01}$,
while the Schwartz kernel of $R''$ takes the form
$$
\rho_{10}^\zeta \rho_{01}^\zeta F'',\quad F''\in\Cinf(\bM^2_0;\pi_R^*\Omega_0);
$$
here $\pi_R^*\Omega_0$ is the lift of the 0-density bundle from the right
factor (a non-vanishing section of which is given by the Riemannian density
$dV_g$). $R(\zeta(k-\zeta))$ is holomorphic, both as a map into the space
of bounded operators on $L^2$ and also into the space of distributions
on $\bM_0^2$, for $\Ree \zeta > k/2$, and extends meromorphically, as a
function with values in the space of distributions on $\bM_0^2$, to
the complex plane when $k$ is even, and to $\Cx \setminus -{\mathbb N}$
when $k$ is odd, with all poles of finite rank. Because it is
constructed using only the symbol calculus, $R'$ is holomorphic in $\lambda$.
Finally, the restrictions of $\rho_{10}^{-\zeta}K_\zeta$ to the left
face $B_{10}$ and of $\rho_{01}^{-\zeta}K_\zeta$ to the right face $B_{01}$
are nonvanishing for all $\zeta$ with $\Ree \zeta > k/2$.
\label{th:ccst}
\end{thm}

\begin{rem}
It will be convenient below to write
$$
K_\zeta = \rho_{10}^\zeta \rho_{01}^\zeta F
$$
where $F$ is smooth on $\bM^2_0$, apart from its conormal singularity
at the lifted diagonal $\Delta \iota_0$.
\end{rem}

For us, the main import of this theorem is its conclusion that $K_\zeta$
simple polyhomogeneous behavior on the space $\bM^2_0$. One of our
ultimate goals here is to find a resolution of the space $X^2$, where $X$ is
the product of two conformally compact manifolds, on which the Schwartz kernel
for the resolvent of its Laplacian is also simple.

In the next sections we shall require uniform weighted $L^2$ estimates for
this resolvent $R(\lambda)$ as $\im\lambda\to\infty$, and we now show how
these follow from the parametrix construction in
(\cite{Mazzeo-Melrose:Meromorphic}). We set $H = -\Delta_g$ here.
\begin{thm} For any $s > 0$ and all $\zeta$ with $\Ree \zeta > k/2 + s$
there is a constant $C_s$ which is independent of $\zeta$ and such that
\begin{equation}
R(\zeta(k-\zeta)): x^s L^2(M,dV_g) \longrightarrow x^s L^2(M,dV_g)
\label{eq:unifwl2o}
\end{equation}
is bounded and satisfies
\begin{equation}
\| R(\zeta(k-\zeta))\|_{\bop(x^s L^2,x^sL^2)} \leq \frac{C_s}{|\im
\zeta(k-\zeta)|}
\label{eq:unifwl2e}
\end{equation}
for all $\zeta$ in this half-plane.
\label{th:unifwl2e}
\end{thm}

\begin{proof}
The proof has a few steps. Setting $\lambda = \zeta(k-\zeta)$, the
parametrix $P(\lambda)$ for $H - \lambda$
constructed in \cite{Mazzeo-Melrose:Meromorphic} satisfies
$$
P(\lambda)(H-\lambda)=\Id+E(\lambda),\quad
(H-\lambda)P(\lambda)=\Id+F(\lambda),
$$
where $E(\lambda)$ and $F(\lambda)$ are residual, and is related to
$R(\lambda)$ by the formula
\begin{equation}\label{eq:parametrix-identity}
R(\lambda)=P(\lambda)-E(\lambda)P(\lambda)+E(\lambda)R(\lambda)F(\lambda).
\end{equation}
The desired uniform estimate for $R(\lambda)$ then follows from
the standard uniform $L^2$ estimate
$$
\|R(\lambda)\|_{\bop(L^2,L^2)}
\leq \frac{1}{|\im\lambda|}, \qquad\lambda\in\Cx\setminus\spec(H),
$$
and from appropriate uniform estimates for $P(\lambda)$, $E(\lambda)$ and
$F(\lambda)$ in weighted spaces which we now derive.

The main work involves demonstrating a uniform estimate in weighted
spaces for the resolvent of the Laplacian in hyperbolic space, and so we
turn to this first. Uniform boundedness of
\[
R(\zeta(k-\zeta)): x^s L^2(\HH^{k+1}; dV) \longrightarrow
x^s L^2(\HH^{k+1}; dV)
\]
is equivalent to the uniform boundedness of the conjugated
operator $x^{-s}R(\zeta(k-\zeta))x^s$ on $L^2(\HH^{k+1}; dV)$.
The Schwartz kernel of this conjugated operator is
\[
K^s_\zeta(z,z') = K^0_\zeta(z,z')(x'/x)^s = K^0_\zeta(z,z')(\rho_{01}/
\rho_{10})^s,
\]
where $K^0_\zeta(z,z')$ is the Schwartz kernel of the resolvent
of the Laplacian in $\HHk$.
Recall from \cite{Mazzeo-Melrose:Meromorphic}
that there is an explicit formula for $K^0_\zeta$:
\begin{eqnarray}
K^0_\zeta(\delta) =
\label{eq:skhs}  \\
c_k\left(\frac{1}{\sinh \delta}
\frac{\del\,}{\del \delta}\right)^{\frac{k-2}{2}}
\left(\frac{1}{\sinh \delta} e^{-(\zeta - k/2)\delta}\right),
&\quad &k\ \mbox{even}, \nonumber \\
c_k\int_0^\infty e^{-(\zeta - k/2)\omega}
(\cosh \omega - \cosh \delta)_+^{-k/2}\, d\omega,  &\quad &k\ \mbox{odd}.
\nonumber
\end{eqnarray}

Here we use that the resolvent is a point-pair invariant, so $R(\lambda;z,z')$
only involves the Riemannian (hyperbolic) distance between $z$ and $z'$,
and equivalently, may be written in terms of the elementary point-pair
invariant $\delta(z,z')$ defined by
\[
\cosh \delta(z,z') = 1 + \frac{|z-z'|^2}{2xx'}.
\]

There are two regions in $(\HH^{k+1})^2_0$ where $K^s_\zeta$ behaves
slightly differently as $\zeta\to\infty$. These are
one near the diagonal, say $\delta(z,z')\leq 1$, and one away from the
diagonal, say $\delta(z,z')\geq 1$. We use a partition of unity
to divide the Schwartz kernel into two parts. It suffices to
prove $L^2$ boundedness of each piece separately. In the former region,
$\delta(z,z')\leq 1$, the factor
$(\rho_{01}/\rho_{10})^s$ is bounded, independent of $\zeta$, hence
the uniform
$L^2$ boundedness of $K^s_\zeta$, and thus \eqref{eq:unifwl2e} can be
proved just as the $L^2$ boundedness of $K^0_\zeta$, directly from
properties of the Schwartz kernel. More explicitly,
the result in this region may be deduced either by invoking the standard 
extension of the symbol calculus with spectral parameter,
or else simply by direct calculation. On the other hand, in
$\delta(z,z')\geq 1$, $K^s_\zeta$ decays exponentially as $\zeta\to\infty$,
so the additional factor $(\rho_{01}/\rho_{10})^s$, which takes the form
$e^{s(\delta(p,z)-\delta(p,z'))}$, where $p$ is some fixed point in
the interior of $\HH^{k+1}$, does not make much difference: the
Schwartz lemma shows the desired bound \eqref{eq:unifwl2e} (and in fact
better bounds for this piece!).

To proceed, we recall that $P(\lambda)$ is constructed in stages, and
as a sum of two terms, $P(\lambda) = P_1(\lambda) + P_2(\lambda)$.
Here $P_1(\lambda)$ is in the small calculus, and is
a slight modification of the operator $P_1'(\lambda)$ obtained from the
use of the symbol calculus to solve away the conormal singularity along
$\Delta \iota_0$. By standard methods, such a $P_1'(\lambda)$ may be
constructed so as to depend holomorphically on $\zeta\in \Cx$, and to have
uniformly bounded norm on $L^2$. It also acts on $x^s L^2$ for each $s$,
with norm depending
on $s$, but not $\zeta$. Next, $P_2(\lambda)$ is obtained after solving
away the first $\ell$ terms of the Taylor series expansion of the
error term $F_1(\lambda) = I - (H-\lambda)P'_1(\lambda)$,
$F_1(\lambda) \in \Psi_0^{-\infty}(M)$, where $\ell$ is sufficiently large
(greater than $s$). This approximate solution of $(H-\lambda)
P_2(\lambda) = F_1(\lambda)$ is found as follows: first solve the normal
equation
\[
(N(H) - \lambda)N(P_2(\lambda)) = N(F_1(\lambda)),
\]
which is the restriction of the equation to the front face $B_{11}$.
But $N(H)$ is naturally identified with the Laplacian on $\HH^{k+1}$,
and so we may apply the result of the discussion above to choose
$P_2'(\lambda)$ with this property, which satisfies the estimate
(\ref{eq:unifwl2e}), and such that $F_2'(\lambda) = F_1(\lambda)
- (H-\lambda)P_1'(\lambda)$ vanishes to first order along $B_{11}$.
Thus
\[
(H-\lambda)(P_1'(\lambda) + P_2'(\lambda)) = I - F_2'(\lambda),
\qquad F_2'(\lambda) \in \Psi^{-\infty,\zeta+1,\zeta}_0(M).
\]
Applying the first $\ell$ terms of the Neumann series for $(I -
F_2'(\lambda))^{-1}$,
\[
I + S^{(\ell)} = I + F_2'(\lambda) + \ldots + (F_2'(\lambda))^{\ell-1}
\]
on the right on both sides of this equation yields
\[
(H - \lambda)(P_1(\lambda) + P_2(\lambda)) = I - F_2(\lambda),
\]
where
\[
P_j(\lambda) = P_j'(\lambda)(I + S^{(\ell)}(\lambda)), \ \  j = 1,2,
\quad
\mbox{and} \qquad \ \ F_2(\lambda) \in \Psi_0^{-\infty,\zeta + \ell,\zeta}(M).
\]
Clearly
\[
\|F_2(\zeta(k-\zeta))\|_{\bop(x^s L^2,L^2)} \leq \frac{C}{|\Imm
\zeta(k-\zeta)|},
\]
with $C$ independent of $\zeta$. We also obtain a remainder term
$E(\lambda)$ from applying $H-\lambda$ on the right of $P(\lambda)
= P_1(\lambda) + P_2(\lambda)$, and this satisfies uniform bounds
(from $L^2$ to $x^sL^2$) of the same form.

These facts taken all together finish the proof.
\end{proof}

\section{Analytic continuation of the resolvent of $H$}
In this section we shall show how to use the integral formula
(\ref{eq:resform}) to obtain an analytic continuation for the resolvent
$R_X(\mu) = (H-\mu)^{-1}$ on $X = M_1 \times M_2$ past the continuous
spectrum of $H = -\Delta_g = H_1 + H_2$. For this continuation we can either
regard $R_X(\mu)$ as an analytic family of bounded operators between weighted
$L^2$ spaces, or else view its Schwartz kernel as an analytic function with
values in an appropriate space of distributions. The key ingredient here is
the existence of similar analytic continuations for the resolvents $R_j(\mu)
= (H_j - \mu)^{-1}$. Although this continuation result holds in considerably
greater generality than just for products of conformally compact spaces,
we shall focus exclusively on this case for the sake of being specific.

We first set up some notation. Let $\dim M_j=k_j+1$. Then it follows
from Section~\ref{sec:conf-compact}, cf. also
\cite{Mazzeo-Melrose:Meromorphic} that the spectrum of $H_j$ decomposes
into the union of a band of continuous spectrum $ [k_j^2/4,+\infty)$,
as well as possibly a finite number of $L^2$ eigenvalues $\lambda_{ji}$,
$i = 1, \ldots, N_j$, $j = 1,2$, in $(0, k_j^2/4)$, with the corresponding
finite rank eigenprojections $\Pi_{ji}$. Next, each $R_j(\mu_j)$ continues
meromorphically from the resolvent set to the Riemann surface $\Sigma$ for
$\sqrt{\mu_j-k_j^2/4}$, which we think of as two copies of $\Cx$ attached
in the usual way along a cut ${\mathcal C}_j$ extending from $k_j^2/4$ to
$\infty$; the resolvent set is identified with the subset of $\Sigma$
where $\im\sqrt{\mu_j-k_j^2/4}<0$. Usually ${\mathcal C}_j$ is taken to
be the ray along the positive real axis, but it will be convenient to choose
it differently later. As already noted, this continuation of $R_j$ is
either as a map into an appropriate spaces of distributions, or else
for any given $s>0$, in the region $\im\sqrt{\mu_j-k_j^2/4}<s$, as bounded
operators from $x^sL^2$ to $x^{-s}L^2$. In any case, its poles are of
finite rank, and we denote them by $\tilde\lambda_{ji}$, with corresponding
finite rank residues $\tilde\Pi_{ji}$ (the tildes are meant to distinguish
these from the eigendata for $H_j$). Note, however, that these poles of
$R_j$ in the nonphysical part of $\Sigma$ need not be simple.

\begin{figure}[ht]
\begin{center}
\mbox{\epsfig{file=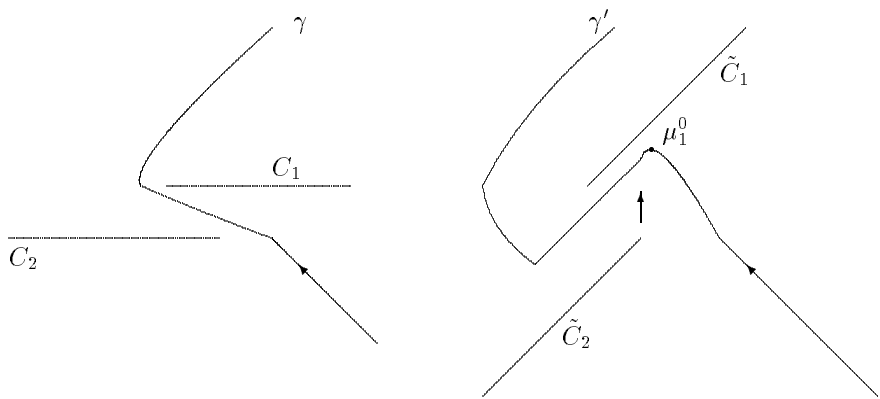}}
\end{center}
\caption{Contour of integration to describe the analytic continuation
of $R(\mu)$. In the second picture, the cuts ${\mathcal C}_1$ and
${\mathcal C}_2$, corresponding to the spectra of $H_1$ and $H_2$,
respectively, have been rotated to the cuts $\tilde{\mathcal C}_1$ and
$\tilde{\mathcal C}_2$. Thus near the point $\mu^0_1$ labeled here,
$R_1(\mu^0_1)$ is {\em not} a bounded operator on $L^2(M_1)$. Now $\mu$,
and hence $\tilde{\mathcal C}_2$ as well, can be moved upwards in the
direction of the arrow.}
\label{fig:contour3}
\end{figure}

We assume until near the end of the argument that neither $H_j$ has any
$L^2$ eigenvalues. Define $k$ by
\[
\frac{k^2}{4} = \frac{k_1^2}{4} + \frac{k_2^2}{4}.
\]
Then $R(\mu)$ is analytic in $\Cx \setminus [k^2/4,\infty)$, and
we wish to show that it continues analytically past the continuous
spectrum. We do this by deforming the contour of integration $\gamma$
in (\ref{eq:resform}) in the following manner. First fix $\mu$ in
the resolvent set, so that $R_1(\mu_1)R_2(\mu - \mu_1)$ is defined
and analytic for $\mu_1 \notin {\mathcal C}_1 \cup (\mu - {\mathcal C}_2)$,
i.e. outside of two horizontal rays, one extending from $k_1^2/4$ to
the right and the other from $\mu - k_2^2/4$ to the left. Next, rotate
these cuts to rays $\tilde{\mathcal C}_j$ by pivoting them by some
angle $\alpha$ counterclockwise around their endpoints. Thus we have
`exposed' two sectors of angle $\alpha$ from the nonphysical portion
of $\Sigma$, and at the same time concealed an equal portion of the
physical part of $\Sigma$. Now it is possible to deform $\gamma$ to
a new contour $\gamma'$ which lies partly in the newly uncovered
sectors in $\Cx\setminus (\tilde{\mathcal C}_1\cup\tilde{\mathcal C}_2)$,
as in Figure~\ref{fig:contour3}. Finally, $\mu$ can then be moved into
the nonphysical region.

We make this a bit more explicit. Suppose that $\mu_0$ is such that
$\mu^0-k_1^2/4 \neq \tilde{\lambda}_{1i}$ and $\mu^0-k_2^2/4 \neq
\tilde{\lambda}_{2i'}$ for any $i,i'$, and also $\mu^0 \neq
\tilde{\lambda}_{1i} + \tilde{\lambda}_{2i'}$. Rotate the cuts ${\mathcal C}_j$
by an angle $\alpha > \arg(\mu^0-k_j^2/4)$, $j = 1,2$. Now deform the
contour $\gamma$ past $\mu^0-k_2^2/4$ to a new contour $\gamma'$, such
that this point is now below $\gamma'$. During the deformation, $\gamma$
may pass through a finite number of poles of $R_1$, yielding residues
$\tilde\Pi_{1i}\otimes R_2(\mu-\tilde\lambda_{1i})$. So long as $\mu$ is
(sufficiently) far away from ${\mathcal C}_2$, $\mu-\mu_1$ is in the
resolvent set of $H_2$, hence no poles of $R_2$ are encountered during
this deformation. Thus we have a new representation
$$
R(\mu)=\frac{1}{2\pi i}\int_{\gamma'}R_1(\mu_1)R_2(\mu-\mu_1)\,d\mu_1
-\sum_{\text{finite}} \tilde\Pi_{1i}\otimes
R_2(\mu-\tilde\lambda_{1i}).
$$
Since $\mu$ is in the resolvent set, the left hand side is still a bounded
operator on $L^2$, although on the right hand side, one factor of the
integrand $R_1(\mu_1)$ is not bounded on $L^2$ along the entire contour.

At this point we merely have a new representation of an operator we
already know exists, namely $R(\mu)$ for $\mu$ in the resolvent set,
in terms of operators which are only bounded between weighted spaces.
However, fixing $\gamma'$, we can now let $\mu$ vary arbitrarily below
this curve, and hence across the spectrum $[k^2/4,\infty)$, as long as the
ramification point $k_2^2/4$ of $R_2$ does not lie on $\gamma'$.
When $\mu-\tilde\lambda_{2i}$ crosses $\gamma'$, the residue term
$-R_1(\mu-\tilde\lambda_{2i})\otimes\tilde\Pi_{1i}$ is produced. This
yields
\[
R(\mu)=\frac{1}{2\pi i}\int_{\gamma'}R_1(\mu_1)R_2(\mu-\mu_1)\,d\mu_1
\]
\[
-\sum_{\text{finite}} \tilde\Pi_{1i}\otimes
R_2(\mu-\tilde\lambda_{1i})
-\sum_{\text{finite}} R_1(\mu-\tilde\lambda_{2i})\otimes\tilde\Pi_{2i};
\]
Each of the residue terms here continue meromorphically in $\mu$ past the
spectrum, with ramification points at $\tilde\lambda_{1i}+k_2^2/4$
(since $R_2$ is ramified at $k_2^2/4$) and similarly, with the indices $1$
and $2$ interchanged. These terms also have poles at $\lambdat_{1i}+
\lambdat_{2j}$, with finite rank residues.

We have proved
\begin{thm} The resolvent $R(\mu)$ for $H$ on $X = M_1 \times M_2$
extends across the cut $[k^2/4,\infty)$ to a meromorphic function
(with values in an appropriate space of distributions) on a Riemann
surface ${\mathcal T}$ ramified at $\tilde\lambda_{1i}+k_2^2/4$ and
$\tilde\lambda_{2i}+k_1^2/4$ (these points are known as Regge poles),
and with finite rank poles at $\lambdat_{1i}+\lambdat_{2i}$.
\end{thm}

When either $H_1$ or $H_2$ has $L^2$ eigenvalues, then the only difference
is that $\gamma'$ must cross these as well. Thus, the analytic continuation
can be written as above, but now we must include sums over the $\lambda_{ji}$
too, and there are new ramification points at $\lambda_{1i} + k_2^2/4$
and $\lambda_{2i} + k_1^2/4$.

\section{Compactification constructions}
We now turn to the other major theme of this paper, which is the various
ways one might compactify the conformally compact spaces $M_j$ or
their product $X = M_1 \times M_2$. The general problem of finding
good compactifications of Riemannian manifolds, and in particular of
(locally) symmetric spaces, has been an area of active research. We refer to
\cite{Ji:Satake} and \cite{Guivarch-Ji-Taylor:Compactifications} for more
discussion of this in the symmetric space setting, and to
\cite{Anderson-Schoen:Positive} and \cite{Schoen-Yau:DG} and
\cite{Freire:Martin} for some beautiful results in some general
geometric settings.

There are two different compactication constructions we shall discuss
here, the geodesic compactification (sometimes also called the conic
compactification), as well as the Martin compactification. Each has a key
role in understanding different aspects of the global geometry and function
theory of a space. We define these now in turn.

When $Z$ is a complete, simply connected manifold of nonpositive curvature
(i.e.a Cartan-Hadamard manifold), then the
geodesic compactification $\bZ$ is obtained by adjoining to $Z$ an
ideal boundary $\del \bZ$, points of which are equivalence classes of geodesic
rays. Two geodesic rays $\gamma(t)$, $\tilde{\gamma}(t)$, $t \geq 0$, are
said to be equivalent if $d(\gamma(t),\tilde{\gamma}(t))$ remains bounded as
$t \to \infty$. Thus in $\RR^n$ any two parallel lines are identified, as are
any two geodesics in the ball model of hyperbolic space $\HH^n$ which converge
to the same point on $S^{n-1}$. If $q \in \del \bZ$, then a neighborhood system
at $q$ is given by sets of the form $(Z \setminus B_R(p')) \cap
\exp_p( [T,\infty) \times {\mathcal U})$, where ${\mathcal U}$ is an
open set in the unit sphere $S_pZ$ in $T_pZ$. A point on this set is
on some geodesic $\gamma(t)$ emanating from $p$ with $\gamma'(0) \in
{\mathcal U}$ and $t \geq T$ as well as in the exterior of the ball
$B_R(p')$. Thus whenever $Z$ is Cartan-Hadamard, then $\bZ$ is
homeomorphic to a closed ball $\overline{B^n}$, and its boundary
$\del \bZ$ is identified with the unit tangent sphere $S_p Z$ for
any $p \in Z$. There are fairly obvious modifications of this construction
when $Z$ is not necessarily simply connected, or only has nonpositive
curvature outside a compact set, and then the structure of $\bZ$ is more
complicated, but only on account of the topology of $\mbox{int}\, Z$.

It is of substantial interest to understand when the compactification
$\bZ$ carries more structure than its initial definition as a
topological space. In particular, we would like to determine when
$\bZ$ is naturally defined as a smooth manifold with boundary,
or with corners. For example, When $Z = \RR^n$ or $\HH^n$, then as we
have already stated $\bZ \approx \overline{B^n}$, but obviously
in these two examples, $\bZ$ is `really' a smooth closed ball.
An interesting and underappreciated feature of this construction is that
although it is possible to identify each of these compactifications with
a ball, the smooth structures at the boundary are different, as we
now explain.

When considering whether $\bZ$ admits a natural smooth structure, the first
key point is the regularity of the transition maps, as we now describe.
Remaining within the context of Cartan-Hadamard manifolds for simplicity, the
transition maps are the homeomorphisms between the unit tangent spheres at
any two points $p,p' \in Z$,
\[
S_p Z \longrightarrow \del \bZ \longleftarrow S_{p'}Z
\]
provided through the geodesic spray from these two points. When $Z$ has
curvatures bounded between two negative constants $-a^2$ and $-b^2$, then
from \cite{Anderson-Schoen:Positive} these transition maps are of
H\"older class ${\mathcal C}^{0,\alpha}$, $\alpha = a/b$, and accordingly,
in this generality, $\bZ$ has only a H\"older structure. However, for either
of the examples above these transition maps are smooth; this is also true
for general conformally compact manifolds for suitable localized versions of
these transition maps (see \cite{Mazzeo:Hodge}). The other key point concerns
the choice of a class of defining functions for $\del \bZ$ (or, if $\bZ$ is
to be regarded as a smooth manifold with corners, then for subsets of it
which are identified with the various boundary hypersurfaces). Smoothness of
the transition maps and choice of defining functions together determine the
${\mathcal C}^\infty$ structure on $\bZ$.

There is considerable
flexibility in choosing (equivalence classes of) defining functions. For
$\RR^n$ it is most natural to use the radial compactification with defining
function for $\del \overline{\RR^n}$, $\rho_1 =1/|z|=1/\mbox{dist}\,(z,0)$,
$z\in\RR^n$, the inverse of the polar distance variable. On the other hand,
the Poincar\'e ball model for $\HH^n$ suggests that the more natural choice
now is $\rho_2 = \exp(-\mbox{dist}\,(z,0))$, $0$ any fixed point in $\HH^n$.
These two defining functions for $\del \overline{B^n}$ are quite different,
since $\rho_1 = -1/\log \rho_2$. We say that $(\overline{B^n},\rho_1)$ is
the `log blow-up' of $(\overline{B^n},\rho_2)$. For later reference, we
shall refer to defining functions defined as the reciprocals of the Riemannian
distance function or its exponential as being of polynomial or exponential
type, respectively.

The difference between these polynomial and exponential type defining
functions appears most clearly in the spaces of functions with
polyhomogeneous behaviour at the boundary, since this involves expansions
in a discrete set of complex powers of the defining function $\rho$ along
with nonnegative integer powers of its logarithm, $\log \rho$. In
particular, functions which are polyhomogeneous with respect to
exponential-type defining functions, with all exponents having positive
real part, are residual, i.e. rapidly decreasing, with respect to
the log blow-up structure. The guiding principle for which defining
functions to choose as primary is determined by the asymptotic behaviour
of eigenfunctions of the Laplacian; one prefers eigenfunctions not to be
automatically residual!  The choices made above for $\overline{\RR^n}$
and $\overline{\HH^n}$, respectively, are vindicated by the analysis
of manifolds with `scattering metrics', see \cite{RBMGeo}, and
with conformally compact metrics, see \cite{Mazzeo-Melrose:Meromorphic} and
\cite{Mazzeo:Hodge}.

As indicated above, the geodesic compactification of the conformally compact
manifold $(M,g)$ is just $\bM$ with its usual smooth structure. It is
obviously of interest to identify the geodesic compactification for the
product of two such manifolds, $X = M_1 \times M_2$. In fact,
$\del \bX$ is the simplicial join of $\pa \bM_1$ and $\pa \bM_2$. More
specifically, it is obtained from
$\del\bM_1\times\del\bM_2\times[0,\pi/2]_\theta$
by collapsing each $\del\bM_1 \times \{q_2\}\times \{0\}$ and
$\{q_1\} \times \del\bM_2 \times \{\pi/2\}$ to a point.
Alternatively, it can be
described as the manifold with corners obtained by blowing up the corner
$\del \bM_1 \times \del \bM_2$ in the product $\bM_1 \times \bM_2$, and
then blowing down the $\bM_2$ fibres in $\del\bM_1\times\bM_2$ and the
$\bM_1$ fibres in $\bM_1\times\del\bM_2$ in the original faces.
There are three boundary hypersurfaces, $\del \bM_1 \times M_2$ and
$M_1 \times \del \bM_2$, corresponding to limits of geodesics of the form
$(\gamma_1(t),p_2)$ and $(p_1,\gamma_2(t))$, respectively, and also the new
face $\del \bM_1 \times \del \bM_2 \times (0,\pi/2)$ corresponding to
limits of geodesics $(\gamma_1(\alpha t), \gamma_2(\beta t))$, $\alpha^2 +
\beta^2 = 1$. It is less clear what defining functions to use for
these faces, though since the `new face' covering the corner comes
from the two-dimensional flats, i.e. products of geodesics in each
factor, it is reasonable that we should use a polynomial type defining
function here. In fact, using the usual coordinates on each factor,
we shall use
\[
\rho_j = -1/\log x_j, \quad j = 1,2, \qquad r = |(\rho_1,\rho_2)|.
\]
In other words, the smooth structure on $\bX$ is the one induced from the
normal blow-up at the corner of the log blow-ups of the two factors $\bM_1$
and $\bM_2$. The reason we use the log blow-up at all faces is that
as we shall prove in the next section, the asymptotic behaviour of
the resolvent involves expansions in powers of these particular
defining functions.

The other compactification we consider here is due to Martin
\cite{Martin:Minimal}, and uses the function theory of $\Lap_g$ to associate
a set of ideal boundary points $\del Z$ to $Z$. Notably, it may be carried
out in great generality for pairs $(Z,H)$ where $H$ is a semibounded
self-adjoint elliptic operator on a space $Z$, though we shall always
assume here that $H$ is the Laplacian. Let $\lambda_0 = \inf\spec(H)$.
Then for every $\lambda \in \RR \setminus (\lambda_0,\infty)$, there
is a compactification $\bZ_M(\lambda)$. Actually, it follows from
the construction that $\bZ_M(\lambda)$ is identified with $\bZ_M(\lambda')$
for any two numbers $\lambda, \lambda' < \lambda_0$, so one needs to
consider only $\bZ_M(\lambda_0)$ and any other $\bZ_M(\lambda)$.

Fix $\lambda < \lambda_0$. By \cite{Sullivan:Related}, the set of solutions
$u$ to the equation $(H - \lambda)u = 0$ which remain everywhere positive
is nonempty; we denote this $\RR^+$-invariant set by
${\mathcal P}_+(\lambda)$. The structure of this positive cone is encoded
in the slice ${\mathcal P}_+^p(\lambda) = \{u \in {\mathcal P}_+(\lambda):
u(p) = 1$ for any fixed $p \in Z$. It follows readily from the Harnack
inequality and elliptic estimates that for any sequence $u_j \in
{\mathcal P}_+^p(\lambda)$ there is a subsequence $u_{j'}$ converging to a
point in this space, and so ${\mathcal P}_+^p(\lambda)$ is compact.
It is also obviously convex. Let ${\mathcal E}$ denote the set of its
extreme points. Then for any $u \in {\mathcal P}_+^p(\lambda)$, the
Krein-Milman theorem gives a measure $dm_u(e)$ supported on ${\mathcal E}$
such that $u = \int_{\mathcal E} dm_u(e)$. This is the generalized
Poisson representation theorem! The set ${\mathcal E}$, or ${\mathcal E}(Z,
\lambda)$ is called the minimal Martin boundary of $Z$. If $u \in
{\mathcal E}$, then whenever $v \in {\mathcal P}_+^p(\lambda)$ and $v \leq u$
then $v = c u$ for some constant $0 < c \leq 1$, and this justifies
the moniker `minimal'.

If $Z$ were to naturally embed in ${\mathcal P}_+^p(\lambda)$, then the
closure of its image would be an obvious way to compactify it. Unfortunately,
this is not the case, but instead we consider the resolvent kernel
$R(\lambda;z,w)$. This is a solution of $(H_z-\lambda)R(\lambda;z,w) = 0$
when $z \neq w$, but is singular at $z = w$ and in addition, $R(\lambda;
p,w) \neq 1$. Thus we define
\[
u_w(z) = R(\lambda,z,w)/R(\lambda,p,w),
\]
so that $u_w(p) = 1$ and $u_w$ is a regular solution of $(H - \lambda)u_w
= 0$ on $X \setminus \{w\}$. Now let $w_j$ be any sequence of points in $Z$
which leaves any compact set. Then some subsequence $u_{w_j'}$ converges to
an element of ${\mathcal P}_+^p(\lambda)$. The (full) Martin boundary is the
set of equivalence classes of these sequences, or equivalently, is the set
of all possible functions $u$ obtained as limits in this fashion. We label
the different boundary points by $q \in \del_M Z(\lambda)$, and write
$u_q(z)$ for the corresponding limiting solution. The Martin compactification
$\bZ_M$ (or $\bZ_M(\lambda)$) is the union of $Z$ and $\del_M Z(\lambda)$.
There is a metric on the set of functions $u_w(z)$, $w \in \bZ_M$, given by
\[
d_p(u_w,u_{w'}) = \int_{B_1(p)} |u_w(z) - u_{w'}(z)|\, dV_g.
\]
Thus $\bZ_M$ not only a topological space, but a metric space.

This definition must be modified when $\lambda = \lambda_0$ and there is an
$L^2$ eigenfunction $u_0$ with eigenvalue $\lambda_0$, since then the
resolvent kernel $R(\lambda_0,z,w)$ does not exist. In this case $u_0 > 0$,
and in fact $\RR^+ \cdot u_0  = {\mathcal P}_+^p(\lambda_0)$, i.e.
the only positive solutions are positive multiples of $u_0$. Hence it is
consistent to let $\bZ_M(\lambda_0)$ be the one point compactification of $Z$
then. Otherwise, if $\lambda_0$ is not in the point spectrum, then the
definition is the same as before.

If $Z = \RR^n$, then $\lambda_0 = 0$ and it is well-known that $\del_M
\RR^n(0)$ is a single point, so that $(\overline{\RR^n})_M(0)$ is the one-point
compactification $S^n$. On the other hand, for $\lambda < 0$, the extreme
points of ${\mathcal P}_+^0(\lambda)$ are the exponentials $e^{x \cdot \xi}$,
$|\xi|^2 = -\lambda$, and this sphere of radius $\sqrt{-\lambda}$ is the full
Martin boundary, and so $(\overline{\RR^n})_M(\lambda) = \bar{B}^n$. If
$Z = \HH^n$ and $n = k+1$, then $\lambda_0 = k^2/4$ and $\bZ_M(\lambda)$ is
the closed ball $\overline{B^n}$ for all $\lambda \leq \lambda_0$. The
minimal positive eigenfunctions are the ones of the form $x^{k/2 +
\sqrt{k^2/4 - \lambda}}$ for all possible choices of upper half-space
coordinates (i.e. choice of which point on the boundary of the ball
to send to infinity). From \cite{Mazzeo-Melrose:Meromorphic} it follows
that that when $M$ is conformally compact, $\bM_M(\lambda)$ is still equal
to $\overline{M}$. To survey other cases relevant to us in which the Martin
compactification is known, when $Z$ is Cartan-Hadamard with curvatures
pinched between two negative constants $-a^2$ and $-b^2$, then Anderson
and Schoen \cite{Anderson-Schoen:Positive} and Ancona \cite{Ancona:Negatively}
proved that $\bZ_M(\lambda) = \overline{B^n}$ (as metric spaces). There has
been recent significant progress in determining the Martin compactifications
for general symmetric spaces of noncompact type; definitive results are proved
in the recent monograph \cite{Guivarch-Ji-Taylor:Compactifications}, and
this has an extensive bibliography of the literature on these developments.

Of particular relevance to us here is the work of Giulini and Woess
\cite{Giulini-Woess:Martin} where the Martin compactification of the
product of two hyperbolic spaces is determined. Their result is that
now $\bX_M(\lambda)$ is identified with $[\overline{\HH^{n_1}}\times
\overline{\HH^{n_2}};\del \overline{\HH^{n_1}}\times \del
\overline{\HH^{n_2}}]$, the blow-up of the product of the two balls
along the corner; recall that this space appeared as an intermediate
picture in the description of the geodesic compactification. Their proof
uses rather involved
global heat-kernel estimates, and one of the motivations of this
paper was to demonstrate how this result, and the analogous one for
products of conformally compact spaces, may be obtained in a
more straightforward manner using resolvent estimates and stationary
phase.

To conclude this section on compactifications, we state once again that
our primary interest is in obtaining a compactification of the double-space
$X \times X$, where $X = M_1 \times M_2$, which is natural with respect to
the resolvent. More specifically, we wish that at least for $\lambda$ in
the resolvent set, $R_X(\lambda)$ should have at most polyhomogeneous
singularities at the boundary hypersurfaces of this compactification
(apart from its usual diagonal singularity). This compactification is
determined by examining the asymptotics of $R_X(\lambda;z,w)$ as the
points $z=(z_1,z_2)$ and $w=(w_1,w_2)$ diverge in all possible directions.
We begin the study these asymptotics in the next section. The Martin
compactification is, at best, essentially a `slice' of this resolvent
compactification, and in any case is easily determined from this asymptotic
analysis, as we do in the final section. We shall see there that when
either $H_1$ or $H_2$ has eigenvalues below the continuous spectrum, then
the Martin compactification of $X$ is obtained by substantially collapsing
part of the boundary of the geodesic compactification $\bX$, and so loses
a lot of information about the fine structure of the resolvent.
Thus the resolvent compactification is the primary object of interest.

\section{Asymptotics}
In the remainder of this paper we shall give a more detailed description
of the structure of the resolvent $R(\mu)$ on $X$. As a first approach
to this we adopt the more traditional viewpoint and derive the asymptotic
behavior of $R(\mu)f$ when $f\in \dCinf(\Xb)$. (Functions vanishing to all
orders at $\del X$ are a suitable analogue of the space of Schwartz functions.)
More generally, it makes absolutely no difference if we allow $f$ to be the
sum of an element of $\dCinf(\Xb)$ and a distribution of compact support.
In particular, all of the calculations below apply when $f = \delta_p$,
$p\in X=M_1\times M_2$; indeed, this is the basis for our identification
of the Martin boundary of $X$. However, we shall simply assume that $f$ is
Schwartz, and also, in this section, that $\mu$ is in the resolvent set
for $H$.

Recall our convention that when $\mu_j$ is in the resolvent set for $H_j$,
then the imaginary part of $\sqrt{\mu_j - k_j^2/4}$ is negative, and by
the results of last section, $R_j(\mu_j)f$ is then of the form
$$
x_j^{k_j/2+i\sqrt{\mu_j-k_j^2/4}}g_j,\ g_j\in\Cinf(\bM_j).
$$
Also, if $\phi_{ji}$ is an $L^2$ eigenfunction of of $H_j$ with eigenvalue
$\lambda_{ji}$, then
$$
\phi_{ji}=x_j^{k_j/2+i\sqrt{\lambda_{ji}-k_j^2/4}}g_j,\ g_j\in\Cinf(\bM_j);
$$
the Schwartz kernel of $\Pi_{ji}$ is $\phi_{ji}(z_j)\otimes
\overline{\phi_{ji}}(z_j')$ if the eigenvalue is simple, and is a finite
sum of such terms otherwise.

We start with a slightly stronger result concerning the structure of
$R_{12}(\mu_1,\mu_2)=R_1(\mu_1)\otimes R_2(\mu_2)$ applied to
$f\in\dCinf(\Xb)$, considered as a function of both $\mu_1$ and $\mu_2$.
Although this follows directly from the corresponding statement in each
factor if $f=f_1\otimes f_2$, $f_j\in\dCinf(\bM_j)$, it is convenient
to prove the general result directly. Thus, the kernel of $R_{12}(\mu_1,
\mu_2)$ is polyhomogeneous on $(\bM_1)^2_0\times (\bM_2)^2_0$.
Let $\pi_L$, $\pi_R$ be the projections to the left and right factors
of $\bM_1\times\bM_2$. That is, if $\pi_{L,j}$, resp.\ $\pi_{R,j}$,
denote the projection
of $(\bM_j)^2_0$ to its left factor, resp.\ right, factor, then
$\pi_L=\pi_{L,1}\times\pi_{L,2}$, $\pi_R=\pi_{R,1}\times\pi_{R,2}$.
Note that $\pi_L$, $\pi_R$, are b-fibrations.
Then $R_{12}(\mu_1,\mu_2)$ applied to $f$ is given by the push-forward
of $R_{12}(\mu_1,\mu_2)\pi_R^* f$ under the map $\pi_L$. Thus
\begin{equation}
R_{12}(\mu_1,\mu_2)f=x_1^{k_1/2+i\sqrt{\mu_1-k_1^2/4}}
x_2^{k_2/2+i\sqrt{\mu_2-k_2^2/4}}g,\ g\in\Cinf(\bM_1\times\bM_2).
\label{eq:jdres}
\end{equation}
The dependence on $\mu_1$ and $\mu_2$ here is uniform in the strong sense
that the function $g(z_1,z_2,\mu_1,\mu_2)$ appearing in
(\ref{eq:jdres}) satisfies
$$
g\in\Cinf(\bM_1\times\bM_2\times
(\Cx\setminus\spec(H_1))_{\mu_1}\times(\Cx\setminus\spec(H_2))_{\mu_2}),
$$
with natural
extensions corresponding to the meromorphic extension of the $R_j$.
This follows directly from the usual push-forward formula \cite{RBMCalcCon}.
{\em Although we use this result throughout the section, we usually state
arguments for simplicity as if $R_1$ and $R_2$ are applied separately to $f$.}
In one particular case, when analyzing $R(\mu)f$ for $\mu$ real, below
$\spec(H)$, and when either $H_1$ or $H_2$ have $L^2$ eigenvalues,
we need a stronger result, where $f$ is not required to be Schwartz.
We postpone that discussion until it is required, see the arguments preceeding
Theorem~\ref{thm:res-set-asymp}.

The asymptotic behavior of $R(\mu)f(z)$ must be analyzed in three separate
regions: near $\del \bM_1 \times M_2$, near $M_1 \times \del\bM_2$
and near the corner $\del \bM_1 \times \del \bM_2$. Using coordinates
$z = (z_1,z_2)$, $z_j = (x_j,y_j)$, these correspond to $x_1 \to 0$,
$x_2 \geq c > 0$, or $x_1 \geq c > 0$, $x_2 \to 0$ or $x_1,x_2 \to 0$,
respectively.

\subsection{Asymptotics at $M_1\times\pa\bM_2$ and at $\pa\bM_1\times M_2$}
We first describe the uniform behaviour of $R(\mu)f$ on $M_1\times
\bM_2$, i.e. at infinity in $M_2$. This is given as an asymptotic series in
powers of $-1/\log x_2$; we shall think of this later as an asymptotic
expansion on the logarithmic blow-up of $\bM_2$, which simply means that we
change the ${\mathcal C}^\infty$ structure of $\bM_2$ by replacing the
defining function $x_2$ for the boundary by $-1/\log x_2$.

\begin{figure}[ht]
\begin{center}
\mbox{\epsfig{file=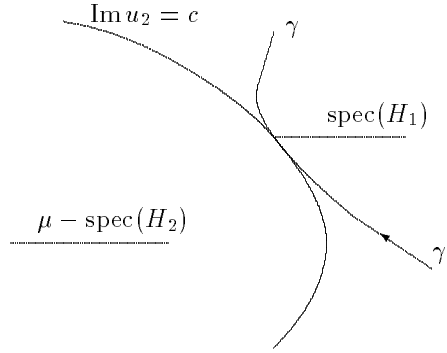}}
\end{center}
\caption{Contour of integration, $\gamma$,
used to describe behavior in $M_1\times
\bM_2$. Here $c=\im\sqrt{\mu-k^2/4}$,
and $u_2=\sqrt{\mu-k_2^2/4-\mu_1}$. \hfill}
\label{fig:contour2}
\end{figure}

Suppose first that $H_1$ has no $L^2$ eigenvalues. We analyze $R(\mu)f$,
using the representation (\ref{eq:resform}) for $R(\mu)$, by shifting
$\gamma$ so that it passes through $\inf\spec(H_1)=k_1^2/4$, and
so that the minimum of $\im\sqrt{\mu-k_2^2/4-\mu_1}$ along this path is
attained at that point, see Figure~\ref{fig:contour2}, and then applying
(complex) stationary phase.

To set things up for stationary phase, first note that $R_2(\mu_2)$ maps
$\dCinf(\bM_2)$ to $x_2^{k_2/2 + i\sqrt{\mu_2-k_2^2/4}}\Cinf(\bM_2)$.
Therefore, the oscillatory factor
$e^{(k_2/2 + i\sqrt{\mu_2-k_2^2/4})\log x_2}$, $\mu_2 = \mu - \mu_1$ appears
in the integrand.
On the other hand, $R_1(\mu_1)$ is {\em not} smooth at $\mu_1=k_1^2/4$, but
instead has an expansion in powers of $\sqrt{\mu_1-k_1^2/4}$. Thus it
decomposes into the sum of odd and even powers, respectively,
$R_1(\mu_1)=R_1^{\text{odd}}(\mu_1)+R_1^{\text{even}}(\mu_1)$, and we
write $R_1^{\text{odd}}(\mu_1)=\sqrt{\mu_1-k_1^2/4}\,G_1(\mu_1-k_1^2/4)$
where $G_1$ is smooth. Because the phase $\sqrt{\mu-\mu_1-k_2^2/4}\,\log x_2$
is not stationary at $k_1^2/4$, the smooth part $R_1^{\text{even}}$
contributes only terms decaying faster than any power of $-1/\log x_2$.
Thus we are left only with
\[
\int_{\gamma} e^{(k_2/2 + i\sqrt{\mu - \mu_1 -k_2^2/4})\log x_2}
\sqrt{\mu_1-k_1^2/4}\,G_1(\mu_1-k_1^2/4) \tilde{g}(\mu_1,z,z')\, d\mu_1,
\]
where $\tilde{g}$ is smooth in all variables. Actually, it suffices to take
this integral only over some compact segment of $\gamma$ containing $k_1^2/4$
and where the decomposition of $R_1$ into even and odd parts is valid;
using the uniform weighted $L^2$ estimates of \S 3, the
integral over the remaining portion of $\gamma$ contributes a term vanishing
like $x_2^{k_2/2 + m}$, $m>0$ can be made as large as wished by
increasing the length of the compact segment,
and is therefore negligible in the
asymptotics. Changing the variable of integration to $\tau =
\sqrt{\mu_1-k_1^2/4}$, so that $\mu_1 = k_1^2/4 + \tau^2$, leads to
\[
\int e^{(k_2/2 + i\sqrt{\mu - k^2/4 - \tau^2})\log x_2}
2\tau^2G(\tau^2)\tilde{g}(k_1^2/4 + \tau^2,z,z')\, d\tau.
\]
The phase function $\phi(\tau) = \sqrt{\mu - k^2/4 - \tau^2}\log x_2$ is
now stationary at $\tau = 0$; furthermore, $\phi''(0) = -1/\sqrt{\mu-k^2/4}$,
and so the further change of variables $\sigma = \tau \sqrt{-\log x_2}/
(\mu-k^2/4)^{1/4}$ reduces this integral to
\[
x_2^{k_2/2 +i\sqrt{\mu-k^2/4}}(\mu-k^2/4)^{3/4}(-1/\log x_2)^{3/2}
\int e^{i\sigma^2/2}\hat{g}(\sigma,z,z')\, d\sigma
\]
where $\hat{g}$ is obtained from $\tilde{g}$ by replacing $\sigma$ by $\tau$.
Applying the stationary phase for phase functions with nonnegative imaginary
part, see \cite[Theorem~7.7.5]{Hor}, we obtain that
\begin{equation}\label{eq:res-side-free-8}
R(\mu)f=
x_2^{k_2/2 + i\sqrt{\mu-k^2/4}}(-1/\log x_2)^{3/2}g,
\end{equation}
where $g$ is $\Cinf$ on the logarithmic blow-up $M_1\times (\bM_2)_{\log}$
of $M_1\times \bM_2$,
i.e. in the variables $(z_1,-1/\log x_2,y_2)$ for $x_1 \geq c > 0$. (Of
course, $g$ is continuous on $M_1\times \bM_2$; it is the lower order terms in
the asymptotics that necessitate the change of the smooth structure.)
In fact,
$$
g|_{M_1\times\pa \bM_2}
= \left. c\, (\mu-k^2/4)^{3/4}x_2^{-k_2/2 -i\sqrt{\mu-k^2/4}}
G_1(0)R_2(\mu-k_1^2/4)f\right|_{M_1\times\pa \bM_2},
$$
where $c$ is a constant arising from the stationary phase lemma.

There is a final simplification of this formula, arising from the fact
that we can identify the operators here slightly more explicitly. First,
according to \cite{Mazzeo-Melrose:Meromorphic}, for each $z_1 \in M_1$,
\[
\left. x_2^{-k_2/2 -i\sqrt{\mu-k^2/4}}R_2(\mu-k_2^2/4)f
\right|_{\pa \bM_2}
\]
\[
= \int_{\pa \bM_2}P_2^t(\mu-k_2^2/4;y_2,z_2')f(z_1,z_2')\, dz_2'
\equiv P_2^t(\mu-k_2^2/4)(f)(z_1,y_2),
\]
where $P_2(\mu-k_2^2/4)$ is the Poisson transform on $M_2$ at the spectral
parameter $\mu-k_2^2/4$, and $P_2^t$ its transpose. Next, by definition,
\[
G_1(0) = \left. \frac{d\,}{d\tau}\right|_{\tau=0}R_1(k_1^2/4 + \tau^2).
\]
We can already see from how it was introduced that this operator arises
from the non-holomorphic part of $R_1$, and so its Schwartz kernel is
smooth on $(\bM_1)^2_0$ because the diagonal singularity of $R_1$,
represented by $R'_1$ in Theorem~\ref{th:ccst}, does not contribute to it.
However, we can see this more directly by differentiating
\[
(H_1 - k_1^2/4-\tau^2)R_1(k_1^2/4 + \tau^2) = I
\]
with respect to $\tau$ and setting $\tau=0$ to get $(H_1 - k_1^2/4)G(0) = 0$,
and so concluding that $G_1(0)$ is smooth by elliptic regularity.
The Schwartz kernel $G_1(0;z_1,z_1')$ is a familiar function on hyperbolic
space: it is a non-square-integrable eigenfunction for $H_1$, with threshold
eigenvalue $k_1^2/4$, invariant under the group of rotations which fix
$z_1'$, and is known as a spherical function. We shall denote it by
$S_1(k_1^2/4;z_1,z_1')$, and regard it as a Schwartz kernel.
Altogether, we have now shown that
\begin{equation}
g|_{M_1\times\pa \bM_2}
= c\, (\mu-k^2/4)^{3/4}S_1(k_1^2/4)P_2^t(\mu-k_2^2/4)f.
\label{eq:x2face}
\end{equation}

When $H_1$ has $L^2$ eigenvalues, then shifting the contour $\gamma$ through
these gives a contribution $-\Pi_{1i}\otimes R_2(\mu-\lambda_{1i})$
from the residues  of the integrand. (For simplicity we assume that the
ramification point $k_1^2/4$ is not a pole of the meromorphic continuation
of $R_1$ to the Riemann surface, though this can obviously be handled too
by the arguments preceeding Theorem~\ref{thm:res-set-asymp}.)
Hence in this case we obtain
\begin{equation}\begin{split}\label{eq:pa-M_2-asymp}
R(\mu)f&=x_2^{k_2/2 + i\sqrt{\mu-k^2/4}}(-1/\log x_2)^{3/2}g\\
&\qquad+\sum_{i=1}^{N_1}x_2^{k_2/2+
i\sqrt{\mu-\lambda_{1i}-k_2^2/4}}g_i,
\qquad \mbox{where}\ \  g,g_i\in\Cinf(M_1\times (\bM_2)_{\log}).
\end{split}\end{equation}
In fact, $g_i = \phi_{1i}\otimes \tilde g_i$, with $\tilde g_i\in\Cinf(
\bM_2)$ and $\phi_{1i}$ an eigenfunction of $H_1$ with eigenvalue
$\lambda_{1i}$; indeed,
$g_i|_{M_1\times \pa\bM_2}=-(\Pi_{1i}\otimes P_2^t(\mu-\lambda_{1i}))f$.

In this region, near $M_1\times\del \bM_2$, only the $L^2$ eigenvalues of
$H_1$ play a role in the asymptotics, but perhaps surprisingly, {\em not}
those of $H_2$. Notice that all terms in \eqref{eq:pa-M_2-asymp} which
come from the $L^2$ eigenvalues of $H_1$ dominate the term coming from the
continuous spectrum since $\im\sqrt{\mu-k^2/4}<\im\sqrt{\mu-\lambda_{1i}-
k_2^2/4}$. In addition, the term corresponding to the lowest eigenvalue
$\lambda_{10}$ of $H_1$ dominates all the other terms. This will be
important later in the determination of the Martin boundary in the presence
of bound states.

The asymptotics of $R(\mu)f$ at the other face $\del \bM_1\times M_2$ is
completely analogous. The same calculations lead to the expression
\begin{equation}\begin{split}\label{eq:pa-M_1-asymp}
R(\mu)f&=x_1^{k_1/2 + i\sqrt{\mu-k^2/4}}(-1/\log x_1)^{3/2}g\\
&\qquad+\sum_{i=1}^{N_1}x_1^{k_1/2+
i\sqrt{\mu-\lambda_{2i}-k_1^2/4}}g_i,
\qquad \mbox{where}\ \  g,g_i\in\Cinf((\bM_1)_{\log}\times M_2).
\end{split}\end{equation}
In addition,
\begin{equation}
g|_{\pa \bM_1\times M_2}
= c\, (\mu-k^2/4)^{3/4}P_1^t(\mu-k_1^2/4)S_2(k_2^2/4)f,
\label{eq:x1face}
\end{equation}
where $P_1^t$ is the transpose of the Poisson transform on $M_1$
and $S_2$ is the `spherical function' at eigenvalue $k_2^2/4$ for
$M_2$.

We note once again that even if both $H_1$ and $H_2$ have $L^2$
eigenvalues, only those of $H_2$ contribute to the asymptotics when
$x_1 \to 0$, $x_2 \geq c > 0$ and only those of $H_1$ contribute to the
asymptotics when $x_2 \to 0$, $x_1 \geq c > 0$. This indicates that the
asymptotics at the corner $\del\bM_1\times \del\bM_2$, where both
$x_1,x_2 \to 0$ must be more complicated, at least in the presence of bound
states, because it intermediates this transition.

\subsection{Asymptotics at $\pa\bM_1\times\pa\bM_2$ in the
absence of $L^2$ eigenvalues}
We now proceed to the analysis of $R(\mu)f$ at this corner.
We have already seen that necessity of logarithmic blow ups at the side faces
of the product. Correspondingly, the asymptotics at the
corner necessitate that we pass to the log blow-ups of both factors
$\bM_j$, which we denote $(\bM_j)_{\log}$. In terms of these, define
\[
\Xt=[(\bM_1)_{\log}\times(\bM_2)_{\log};\pa(\bM_1)_{\log}
\times\pa(\bM_2)_{\log}].
\]
This is the correct space on which to consider asymptotics of $R(\mu)f$, as
we shall now prove.

We let
$$
\rho_j=-1/\log x_j,\quad j=1,2,
$$
be boundary defining functions of $(\bM_j)_{\log}$, and we keep denoting
their pull-back to $\Xt$ with the same notation. Also let $y_1$, $y_2$
denote coordinates on $\pa(\bM_j)_{\log}$, $j=1,2$. Thus, $\Xt$ is the blow-up
of $(\bM_1)_{\log}\times(\bM_2)_{\log}$ at the corner $\rho_1=\rho_2=0$.
Hence valid coordinates in the region where $\rho_1/\rho_2<C$, for some
$C>0$, are given by
$$
\rho_2,\ s=\rho_1/\rho_2,\ y_1,\Mand y_2.
$$
In terms of the original boundary defining functions $x_j$ this means that
in any region where $\log x_2/\log x_1<C$, for some $C>0$, we use the
projective coordinate
\[
s = \frac{-1/\log x_1}{-1/\log x_2} = \frac{\log x_2}{\log x_1}
\]
along the new `front' face of $\Xt$ covering the corner, and $-1/\log x_2$
as a defining function for this face. Thus $s\to 0$ upon approach to the
lift of $\pa \bM_1\times \bM_2$, while $s\to\infty$ on approach to the lift of
$\bM_1\times\pa \bM_2$. We note that a {\em total} boundary defining function
of $\Xt$ is given by
$$
\rho=(\rho_1^{-2}+\rho_2^{-2})^{-1/2},
$$
i.e.\ with the usual Euclidean notation, $r_j=\rho_j^{-1}$, $r=\rho^{-1}$,
$r=\sqrt{r_1^2+r_2^2}$. This explains the appearance of $\rho^{-1}$ in
our results below. In the region $s<C$, $\rho_1=-1/\log x_1$ is another total
boundary defining function; it will be used for the initial local calculation
and then we change to $\rho$ for the global statements. The trivial identity
\[
x_2 = e^{\log x_2} = e^{s \log x_1} = x_1^s
\]
will be used repeatedly.

The integrand $R_1(\mu_1)R_2(\mu-\mu_1)f$ in the contour integral
representation of $R(\mu)f$ has the form
\[
x_1^{k_1/2+i\sqrt{\mu_1-k_1^2/4}}x_2^{k_2/2+i\sqrt{\mu-\mu_1-k_2^2/4}}g,
\]
or equivalently,
\begin{equation}
\label{eq:integrand-5}
\exp\left\{\left[k_1/2+s k_2/2+i\left(\sqrt{\mu_1-k_1^2/4}+
s\sqrt{\mu-\mu_1-k_2^2/4}\right)\right]\log x_1\right\}g,
\end{equation}
where $g$ is $\Cinf$ in $\mu, \mu_1$ and on $\bM_1 \times \bM_2$.
The expressions inside these square roots assume values in $\Cx\setminus
[0,\infty)$, and we assume the square roots have {\em negative} imaginary
parts in this region.

Suppose first that neither $H_1$ nor $H_2$ has $L^2$ eigenvalues. We again
do a stationary phase analysis of the integral of (\ref{eq:integrand-5}),
and for this it is necessary to choose the contour of integration so that,
if we set
\[
F(\mu_1)=F_s(\mu_1) =  \sqrt{\mu_1-k_1^2/4}+s\sqrt{\mu-\mu_1-k_2^2/4},
\]
then the supremum of $\im F(\mu_1)$ along $\gamma$ is as negative as possible.
Both $\mu$ and $s$ are parameters here, and we must choose the contour
differently corresponding to the different points of the front face of $\Xt$,
i.e. the different values of $s$. Since $F$ is analytic outside $[k_1^2/4,
+\infty) \cup(\mu-[k_2^2/4, +\infty))$, its critical points and those of its
imaginary
part are the same. In fact, $F$ has a unique critical point, located at
\begin{equation}\label{eq:crit-point-im}
\mu_1^0=\mu_1^0(s)=\frac{\mu-k_2^2/4+s^2 k_1^2/4}{1+s^{2}},
\end{equation}
and this is a saddle point of the harmonic function
$\im F(\mu_1)$. (We return to this and shall
explain it in greater detail below.) Moreover, this critical point
always lies on the straight line segment connecting the two ramification
points $k_1^2/4$ and $\mu-k_2^2/4$, and tends to the former as $s\to+\infty$
and to the latter as $s\to 0$. Next,
\[
F(\mu_1^0) = \sqrt{(\mu-k^2/4)(1+s^2)}, \qquad F'(\mu_1^0) = 0,
\]
and
\[
F''(\mu_1^0) = -\frac{1}{4}s^{-2}(1+s^2)^{5/2}(\mu-k^2/4)^{-3/2}.
\]
Therefore we may choose an appropriate contour $\gamma$ so that
$\im F(\gamma(t))$ attains a non-degenerate maximum when $\gamma(t) =
\mu_1^0$. The stationary phase lemma may now be applied as before, although
the singular change of variables is no longer required (in $s>0$);
we deduce that the
asymptotics of $R(\mu)f$ have the form
\begin{equation}
\begin{split}\label{eq:corner-asymp-free-8}
&\exp\left(\left(k_1/2+s k_2/2+i\sqrt{(\mu-k^2/4)(1+s^2)}\right)
\log x_1\right)(-1/\log x_1)^{1/2} g\\
&=x_1^{k_1/2}x_2^{k_2/2}\left(\exp\sqrt{(\log x_1)^2+(\log x_2)^2}\right)
^{-i\sqrt{\mu-k^2/4}}(-1/\log x_1)^{1/2}g.
\end{split}
\end{equation}
At the points of the front face of $\Xt$ where $\log x_2/\log x_1 =s$
this coefficient function $g$ is a nonvanishing multiple of
\[
s(1+s^2)^{-5/4}(\mu-k^2/4)^{3/4}g(\mu_1^0(s)) = s\, (1+s^2)^{-5/4}
((\mu-k^2/4)^{-3/2})^{-1/2}
\]
\[
\times x_1^{-k_1/2-i\sqrt{\mu_1^0(s)-k_1^2/4}}
x_2^{-k_2/2-i\sqrt{\mu-\mu_1^0(s)-k_2^2/4}}R_1(\mu_1^0(s))R_2(\mu-\mu_1^0(s))f,
\]
or, finally,
\begin{equation}
c \, s\, (1+s^2)^{-5/4}(\mu-k^2/4)^{3/4}
P_1^t(\mu_1^0(s))P_2^t(\mu-\mu_1^0(s))f.
\label{eq:frontface}
\end{equation}
Note that although this computation is done separately for each value of
$s$, this final expression depends smoothly on $s$ and in fact $g$ is smooth
on $\Xt$ up to the front face.

To be precise, we still need to examine what happens near and at the corner
$s=0$, $-1/\log x_2=0$. However, this is hardly different
from the discussion at $\pa\bM_1\cap M_2$. Indeed, we simply introduce
$\tau=\sqrt{\mu-\mu_1-k_2^2/4}$ as the new smooth variable of integration,
and then stationary phase can be performed uniformly in $s$ down to $s=0$,
with the limiting contour being the same as for the
discussion at $\pa\bM_1\cap M_2$, showing that
$R(\mu)f$ is actually polyhomogeneous. In particular,
note that the expansions (\ref{eq:x2face}), (\ref{eq:x1face}) and
(\ref{eq:frontface}) at the front and side faces of $\Xt$ match up at the
corners, in agreement with the polyhomogeneity.
The one point to note is that, for
example, there is a factor of $(-1/\log x_1)^{1/2}$ in (\ref{eq:frontface}),
whereas (\ref{eq:x1face}) has a factor of $(-1/\log x_1)^{3/2}$.
To reconcile this, observe that (\ref{eq:x1face}) has an extra factor of
$s$, and
\[
s (-1/\log x_1)^{1/2} = \frac{-\log x_2}{-\log x_1}(-1/\log x_1)^{1/2}
= (-\log x_2)(-1/\log x_1)^{3/2},
\]
so in fact the powers match up.

We remark that if the contour $\gamma$ does not go through the critical
point $\mu_1^0(s)$, then it must necessarily contain points where the
integrand is larger; however, at those points the phase function itself,
{\em not just its imaginary part}, will fail to be stationary, and so
stationary phase gives a decay rate $O((-1/\log x)^{\infty})$; this
merely indicates that the contour has not been chosen optimally and a
better result is possible.

Before proceeding to the rather more subtle discussion of asymptotics
at the corner in the presence of bound states, we summarize
our results in their absence.

\begin{prop}\label{prop:no-bd-state-asymp}
Suppose $f\in\dCinf(\bX)$, $\mu\in \Cx\setminus[k^2/4,+\infty)$, and
$H_1$, $H_2$ have no $L^2$ eigenfunctions. Then, with
$\rho_i=-1/\log x_i$, $\rho^{-1}=\sqrt{\rho_1^{-2}+\rho_2^{-2}}$,
$R(\mu)f$
has the following asymptotic expansion on $\Xt$:
\begin{equation}\begin{split}
R(\mu)f&=x_1^{k_1/2}x_2^{k_2/2}\exp (-i\sqrt{\mu-k^2/4}/\rho)h.
\end{split}\end{equation}
Here $h$ polyhomogeneous on $\Xt$, with order
$1/2$ on the front face and $3/2$ on the side faces of $\Xt$.
Moreover, the principal symbol of $h$ is an elliptic $f$-independent
multiple of
\begin{equation}\begin{split}
&P_1^t(\mu_1^0(s))P_2^t(\mu-\mu_1^0(s))f\ \text{on the front face},\\
&
S_1(k_1^2/4)P_2^t(\mu-k_2^2/4)f\ \text{on the lift of}\ \bM_1\times
\pa\bM_2,\\
&P_1^t(\mu-k_1^2/4)S_2(k_2^2/4)f\ \text{on the lift of}\ \pa\bM_1\times \bM_2,
\end{split}\end{equation}
with $\mu_1^0(s)$ given by \eqref{eq:crit-point-im}.
\end{prop}

This discussion already indicates that $s$ plays a dual role, both as
a coordinate on the space $\tilde{X}$ and also as a spectral parameter.
The function $\mu_1^0(s)$ identifies the front face of $\Xt$ with
the line joining the two threshold values $k_1^2/4$ and $\mu-
k_2^2/4$ in the spectral plane. This role becomes even more pronounced
in determining asymptotics at the front face in the presence of bound
states, for then the optimal locus for the contour $\gamma$ determines
which of the residue terms corresponding to the different eigenvalues
must be included in the asymptotics. In addition, for certain real
values of $\mu$, there are a finite number of exceptional values of $s$
at which $\mu_1^0(s)$ equals either $\lambda_{1i}$ or $\mu - \lambda_{2i}$;
then the optimal contour must pass through this pole, and this creates an
additional Legendrian singularity. We explain all of this now.

\subsection{Asymptotics at $\pa\bM_1\times\pa\bM_2$ in the
presence of $L^2$ eigenvalues for $\mu\in\Cx\setminus\Real$}
Thus suppose that $H_1$ or $H_2$ have $L^2$ eigenvalues. As we have
already seen, these correspond to poles of the integrand and the
residues at these poles may contribute to the asymptotics. To see
when this happens, first note that the additional terms coming from
these residues have the form
\begin{equation}\label{eq:corner-asymp-ev-8}
x_2^{k_2/2+
i\sqrt{\mu-\lambda_{1i}-k_2^2/4}}\phi_{1i}\otimes \tilde g_{1i},\ \tilde g_{1i}
\in\Cinf(\bM_2),
\end{equation}
with analogous terms when the roles of $\bM_1$ and $\bM_2$ have been switched.
Here
$$
\phi_{1i}\in x_1^{k_1/2+i\sqrt{\lambda_{1i}-k_1^2/4}}\Cinf(\bM_1).
$$
The first main observation is that the terms \eqref{eq:corner-asymp-ev-8}
dominate the previous one \eqref{eq:corner-asymp-free-8} coming from
the continuous spectrum only when
$$
\im\left(\sqrt{\lambda_{1i}-k_1^2/4}+s\sqrt{\mu-\lambda_{1i}-k_2^2/4}\right)
> \im\left(\sqrt{\mu-k^2/4}\sqrt{1+s^2}\right),
$$
or in other words when
\begin{equation}
\im F(\lambda_{1i}) > \im F(\mu_1^0(s)).
\label{eq:eigineq}
\end{equation}
In all other cases, these terms coming from the bound states are lower
order in the expansion \eqref{eq:corner-asymp-free-8}.

To analyze this, we must examine the geometry of the function $F$ a bit
more closely. Set
\[
\Lambda(\mu) = \Cx \setminus \big( [k_1^2/4,\infty)\cup (\mu - [k_2^2/4,
\infty)\big),
\]
and let $L \subset \Lambda(\mu)$ be the open line segment connecting the
two threshold values $k_1^2/4$ and $\mu-k_2^2/4$; this segment is the
image of the map $(0,\infty) \ni s \to \mu_1^0(s)$. The function $F$
depends implicitly on both $\mu$ and $s$; we shall hold $\mu$ fixed
throughout, but shall study how $F$ changes as $s$ varies.  We have
already noted that $\im F$ has a saddle point at $\mu_1^0(s)$. To see this,
observe first that this is the unique critical point. Next,
\[
\lim_{\mu_1 \to k_1^2/4} \im F(\mu_1) = s \im \sqrt{\mu - k^2/4},
\qquad
\lim_{\mu_1 \to \mu-k_2^2/4} \im F(\mu_1) = \im \sqrt{\mu - k^2/4},
\]
and both of these values are greater than $\im F(\mu_1^0(s)) =
\sqrt{1+s^2}\im \sqrt{\mu-k^2/4}$ since $\im \sqrt{\mu-k^2/4} < 0$.
Also, each of the square roots in the expression for $F$ has imaginary
part vanishing along one of the slits and tending to $-\infty$ as $\mu_1$
moves horizontally away from that slit, and so $\im F$ tends to
$-\infty$ as $|\mu_1| \to \infty$. Thus $\mu_1^0(s)$ is indeed
a saddle point for $\im F$.
We make the following definition.

\begin{Def}\label{Def:Ns-Ps}
For each $s$, let ${\mathcal N}_s$ and
${\mathcal P}_s$ denote the regions in $\Lambda(\mu)$ where
$\im F(\mu_1) < \im F(\mu_1^0(s))$ and $\im F(\mu_1) > \im F(\mu_1^0(s))$,
respectively.
\end{Def}

The union of ${\mathcal N}_s$ and the slits always has one component, while
${\mathcal P}_s$ has two components, each of which has compact closure
intersecting precisely one of the two slit axes and which intersect at
$\mu_1^0(s)$. Label these two components ${\mathcal P}_s^\ell$ and
${\mathcal P}_s^r$, respectively, the $\ell$ and $r$ denoting whether
these regions touch the slits extending toward the left or right.
Finally, recall that the `optimal' contour of integration $\gamma$ is
any contour which remains entirely within the region ${\mathcal N}_s$,
except when it passes through the saddle point.

\begin{figure}[ht]
\begin{center}
\mbox{\epsfig{file=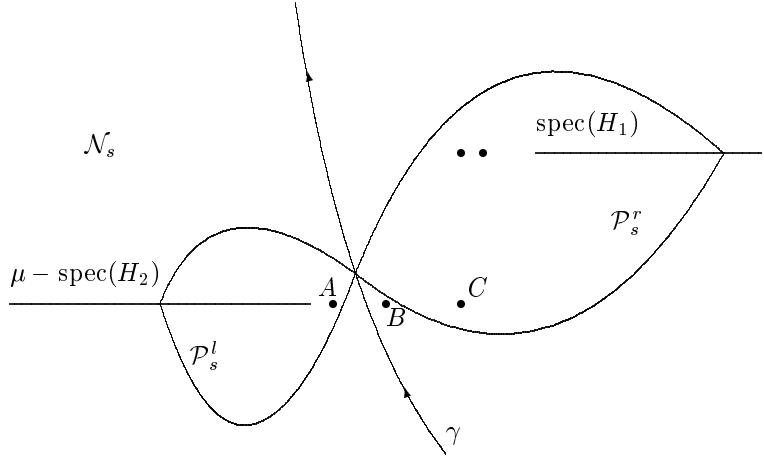}}
\end{center}
\caption{The regions ${\mathcal P}_s^\ell$, ${\mathcal P}_s^r$ and
${\mathcal N}_s$ for small $s>0$. The points $A$, $B$ and $C$ denote
$\mu-\lambda_{23}$, $\mu-\lambda_{22}$ and $\mu-\lambda_{21}$ respectively.
The contour $\gamma$ had to be shifted
through point $C$, but not yet past $A$. The point $B$ can made to lie
on either side of $\gamma$; the residue there gives a smaller asymptotic
term than the critical point $\mu_1(s)$ of $F$.}
\label{fig:contour4}
\end{figure}

Now consider the locations of the poles $\lambda_{1i}$ and $\lambda_{2i}$
relative to the contour $\gamma$ and these regions ${\mathcal N}_s$
and ${\mathcal P}_s$. For any $s$, these two collections of poles are
separated by the initial choice of contour $\gamma$ in \eqref{eq:resform}.
However, in deforming $\gamma$ to a $\gamma'$ which lies entirely
within ${\mathcal N}_s$ it may be necessary to shift past some of
these poles. This is not necessary if all of the $\lambda_{1i}$
lie within ${\mathcal N}_s \cup {\mathcal P}_s^r$ and all of the
$\mu-\lambda_{2i}$ lie within ${\mathcal N}_s \cup {\mathcal P}_s^\ell$.
When $s$ is sufficiently small, then in fact all of the $\lambda_{1i}$
lie within ${\mathcal N}_s \cup {\mathcal P}_s^r$; all of the
$\mu-\lambda_{2i}$ lie in this region too, and none lie in
${\mathcal P}_s^\ell$. On the other hand, when $s \to \infty$,
both the $\lambda_{1i}$ and the $\lambda_{2i}$ lie within
${\mathcal N}_s \cup {\mathcal P}_s^\ell$. As $s$ increases from $0$
to $\infty$, the region ${\mathcal P}_s^\ell$ gradually engulfs each of
the points $\lambda_{1i}$, and so the optimal contour (which must always
lie in ${\mathcal N}_s$) must cross these before this happens. A similar
phenomenon occurs for the $\lambda_{2i}$.

\subsection{Asymptotics at $\pa\bM_1\times\pa\bM_2$ in the
presence of $L^2$ eigenvalues for $\mu\in\Real$}
On the other hand, when $\mu \in \RR$ and there are bound states in the
interval $(\mu-k_2^2/4,k_1^2/4)$, i.e. when
\[
\mu < \min\big\{k_1^2/4 - \max\{\lambda_{1i}\},
k_2^2/4 - \max\{\lambda_{2i}\}\big\},
\]
then this analysis must fail for certain values of $s$, namely those
values where either $\mu_1^0(s) = \lambda_{1i}$ or $\mu - \mu_1^0(s) =
\lambda_{2i'}$ for some $i,i'$. For then an optimal contour would need to
pass directly through a pole. It is clear that there must be some sort
of jump in the behaviour at these points, because for nearby values of $s$
this phenomenon does not occur and the previous analysis may be used to get
the asymptotics; the residue at the pole would be included in these
asymptotics if $s$ varies slightly to one side, but are not included if it
varies slightly to the other. We now prove that $R(\mu)$ has a Legendrian
singularity at this front face of $\Xt$ at those values of $s$ where this
occurs.

To be definite, suppose that $s_0$ satisfies $\mu_1^0(s_0) =\lambda_{1i}$.
Let $\gamma$ be the original contour defining $R(\mu)$, as in
Figure~\ref{fig:contour}, and let $\gamma_s$ be an optimal contour for any
value $s$ near $s_0$. We may arrange that $\gamma_s(0)=\mu_1^0(s)$,
so $\gamma_{s_0}(0)=\lambda_{1i}$.
When $s>s_0$ (closer to the lift of
$\bM_1\times\pa \bM_2$) the contour has crossed past $\lambda_{1i}$, and
the residue term is present, whereas when $s<s_0$ (closer to the lift of
$\pa \bM_1\times \bM_2$) this residue term is not included.

We proceed as follows. First consider what happens when $s\to s_0$ from
the left, i.e. with $s\leq s_0$. Equation \eqref{eq:corner-asymp-free-8}
describes the asymptotics when $s<s_0$, and we are interested in the
the behavior in the limit. Write
$$
R_2(\mu-\mu_1)=R_2(\mu-\lambda_{1i})+(\mu_1-\lambda_{1i})\tilde R_2(\mu-\mu_1),
$$
i.e.\ expand $R_2$ to first order in Taylor series in $\mu_1$ around
$\mu_1^0(s_0)=\lambda_{1i}$. Thus,
\begin{equation}\begin{split}\label{eq:tilde-R_2-def}
\tilde R_2(\mu-\mu_1)=
\frac{R_2(\mu-\mu_1)-R_2(\mu-\lambda_{1i})}{\mu_1-\lambda_{1i}} \\
=-\int_0^1 R_2'(\sigma(\mu-\mu_1)+(1-\sigma)(\mu-\lambda_{1i}))\,d\sigma,
\end{split}
\end{equation}
and in particular, $\lim_{\mu_1\to\lambda_{1i}}\tilde R_2(\mu-\mu_1)=
-R'_2(\mu-\lambda_{1i})$.

Our first remark is that the contour in
$\int_{\gamma_s} R_1(\mu_1)R_2(\mu-\lambda_{1i})f\,d\mu_1$
can be shifted farther from the spectrum of $H_1$ since now the $R_2$ term is
independent of $\mu_1$, so its contribution is negligible. Hence
we only need to analyze
$\int_{\gamma_s} R_1(\mu_1)(\mu_1-\lambda_{1i})\tilde R_2(\mu-\mu_1)\,d\mu_1$.
But
$$
\tilde R_1(\mu_1)\tilde R_2(\mu-\mu_1)
=R_1(\mu_1)(\mu_1-\lambda_{1i})\tilde R_2(\mu-\mu_1)
$$
is holomorphic for $\mu_1$ near $\lambda_{1i}$, so we may move the contour
to go through $\lambda_{1i}$ when $s=s_0$.
Moreover, from the usual parametrix identity
\eqref{eq:parametrix-identity}, using the holomorphy of $\tilde R_1(\mu_1)
=(\mu_1-\lambda_{1i})
R_1(\mu_1)$ as a bounded operator on $L^2$ near $\mu_1=\lambda_{1i}$,
$(\mu_1-\lambda_{1i})R_1(\mu_1)$ has the standard asymptotics. Note
that $\tilde R_1(\lambda_{1i})=\Pi_{1i}$.

Thus,
the holomorphic function $\tilde R_1(\mu_1)R_2(\mu_2)f$ has asymptotics
$$
\tilde R_1(\mu_1)R_2(\mu_2)f=x_1^{k_1/2}x_2^{k_2/2}x_1^{i\sqrt{\mu_1-k_1^2/4}}
x_2^{i\sqrt{\mu_2-k_2^2/4}} a(\mu_1,\mu_2),
$$
where $a(\mu_1,\mu_2)$ is $\Cinf$ on $U_{\mu_1}\times U_{\mu-\mu_2}
\times \bM_1\times\bM_2$, and holomorphic in $\mu_1$ and $\mu_2$,
where $U$ is a neighborhood of $\lambda_{1i}$. Hence
\[
\tilde R_1(\mu_1)\tilde R_2(\mu-\mu_1)f = \frac{\tilde R_1(\mu_1)
R_2(\mu-\mu_1) f -\tilde R_1(\mu_1)R_2(\mu-\lambda_{1i})f}{\mu_1
-\lambda_{1i}}
\]
\[
=x_1^{k_1/2}x_2^{k_1/2} \cdot (\mu_1-\lambda_{1i})^{-1} \cdot
\left(x_1^{i\sqrt{\mu_1-k_1^2/4}}x_2^{i\sqrt{\mu-\mu_1-k_2^2/4}}
a(\mu_1,\mu-\mu_1) \right.
\]
\[ \left. - x_1^{i\sqrt{\mu_1-k_1^2/4}}
x_2^{i\sqrt{\mu-\lambda_{1i}-k_2^2/4}}
a(\mu_1,\mu-\lambda_{1i}) \right)
\]
Now add and subtract $x_1^{i\sqrt{\mu_1-k_1^2/4}}
x_2^{i\sqrt{\mu-\mu_1-k_2^2/4}} a(\mu_1,\mu-\lambda_{1i})$ in the numerator
of this fraction on the right, and note that the term involving the
holomorphic function
$$
\frac{a(\mu_1,\mu-\mu_1)-a(\mu_1,\mu-\lambda_{1i})}{\mu_1-\lambda_{1i}}
$$
has the same asymptotics as if $\lambda_{1i}$ were not a pole; thus we
only need to consider
$$
\int_{\gamma_s}\frac{x_1^{i\sqrt{\mu_1-k_1^2/4}}x_2^{i\sqrt{\mu-\mu_1-k_2^2/4}}
-x_1^{i\sqrt{\mu_1-k_1^2/4}}
x_2^{i\sqrt{\mu-\lambda_{1i}-k_2^2/4}}}{\mu_1-\lambda_{1i}}\,a(\mu_1,
\mu-\lambda_{1i})\,d\mu_1.
$$
We can also write $a(\mu_1,\mu-\lambda_{1i})=a(\lambda_{1i},\mu-
\lambda_{1i})+(\mu_1-\lambda_{1i})a_1(\mu_1,\mu-\lambda_{1i})$, where
$a_1$ is holomorphic in $\mu_1$ and smooth on $\bM_1\times\bM_2$;
the term involving $a_1$ again yields an expression with the
standard asymptotics, as if $\lambda_{1i}$ were not a pole.
Hence it remains to consider
$$
a(\lambda_{1i},\mu-\lambda_{1i})
\int_{\gamma_s}\frac{x_1^{i\sqrt{\mu_1-k_1^2/4}}x_2^{i\sqrt{\mu-\mu_1-k_2^2/4}}
-x_1^{i\sqrt{\mu_1-k_1^2/4}}
x_2^{i\sqrt{\mu-\lambda_{1i}-k_2^2/4}}}{\mu_1-\lambda_{1i}}\,d\mu_1.
$$
This depends on $f$ only via
\begin{equation*}\begin{split}
a(\lambda_{1i},\mu-\lambda_{1i})&=-x_1^{-k_1/2-i\sqrt{\lambda_{1i}-k_1^2/4}}
x_2^{-k_2/2-
i\sqrt{\mu-\lambda_{1i}-k_2^2/4}}(\Pi_{1i}\otimes R_2(\mu-\lambda_{1i}))f\\
&=-x_1^{-k_1/2-i\sqrt{\lambda_{1i}-k_1^2/4}}
(\Pi_{1i}\otimes P_2^t(\mu-\lambda_{1i}))f.
\end{split}\end{equation*}
The factor $x_1^{-i\sqrt{\lambda_{1i}-k_1^2/4}}$ here corresponds
to the precise decay of the eigenfunction $\phi_{1i}$.
The remaining integral is an explicit function which depends only
on $s$, $x_1$ and $x_2$ (and of course $\mu$). Note that in the region $x_2>0$,
it is simply a smooth function of $s$, and $x_2$ (or equivalently,
$s$ and $x_1$, keeping in mind that $s=\log x_2/\log x_1$ is near
$s_0\in(0,+\infty)$. Thus we only need to understand the asymptotics of
this function as $x_1,x_2\to 0$.

Factor out $x_1^{iF(\mu_1)}=x_1^{iF(\mu_1^0(s))}
x_1^{i(F(\mu_1)-F(\mu_1^0(s)))}$ to get
\begin{equation}\label{eq:integral-77}
x_1^{iF(\mu_1^0(s))}\int_{\gamma_s}e^{-i(F(\mu_1)-F(\mu_1^0(s)))/\rho_1}
\frac{1-e^{-is(\sqrt{\mu-\lambda_{1i}-k_2^2/4}-\sqrt{\mu-\mu_1-k_2^2/4})
/\rho_1}}{\mu_1-\lambda_{1i}}\,d\mu_1,
\end{equation}
where we wrote the integrand in terms of $\rho_1=-1/\log x_1\geq 0$ to keep
the signs easier to follow.
Note that $\im(F(\mu_1)-F(\mu_1^0(s)))\leq 0$ along $\gamma_s$.

So far {\em any} contour $\gamma_s$
that stays in ${\mathcal N}_s$ has been suitable
for our calculations. Now we impose an additional condition, namely
that there exists a fixed interval $[-\ep,\ep]$ such that
for $s_0-\delta\leq s\leq s_0$ $\gamma_s$ stays in the region where
$\im \sqrt{\mu-\mu_1-k_2^2/4}\geq \im \sqrt{\mu-\lambda_{1i}-k_2^2/4}$.
Since this
region is one side of a parabola, see Figure~\ref{fig:contour2}, this
condition can be easily satisfied.
This condition ensures that the real part of the exponent in the
numerator in \eqref{eq:integral-77} is non-positive, hence bounded even
as $\rho_1\to 0$. In addition, the fraction in \eqref{eq:integral-77}
is bounded by $C\rho_1^{-1}$ for some $C>0$, since dividing
the denominator by $\rho_1$ yields a bounded function.

The key point is that
$$
F(\mu_1)-F(\mu_1^0(s))=(\mu_1-\mu_1^0(s))^2F_2(\mu_1-\mu_1^0(s),\mu_1^0(s)),
$$
with $F_2$ non-zero at $(0,\mu_1^0(s))$, since $\mu_1^0(s)$ is
a non-degenerate critical point of $F$, while
\begin{equation*}\begin{split}
\sqrt{\mu-\lambda_{1i}-k_2^2/4}-\sqrt{\mu-\mu_1-k_2^2/4}
=(\mu_1-\lambda_{1i})F_1(\mu_1-\lambda_{1i})
\end{split}\end{equation*}
with $F_1$ nonvanishing and holomorphic near $0$.
This implies that the first exponent is a smooth function of
$(\mu_1-\mu_1^0(s))^2/\rho_1$, with
a negative real part of the same order of magnitude as $\rho_1\to 0$,
hence Schwartz as
$(\mu_1-\mu_1^0(s))^2/\rho_1\to \infty$,
and the second one is a smooth function of $(\mu_1-\lambda_{1i})/\rho_1$ with
a negative real part. However, we can only make the Schwartz conclusion here
if the real part of the same order of magnitude as
$(\mu_1-\lambda_{1i})/\rho_1$. Since $\gamma_s$ is on one side of a parabola,
this is impossible to accomplish if $\gamma_s$ is to be smooth (for then
it is tangent to the parabola). Hence, we need to break up $\gamma_s$ into
two integrals which will yield boundary terms.
These boundary terms cancel when $s\neq s_0$.

The geometry is as follows. Start with the space
$Z_0=[0,1)_{\rho_1}\times(s_0-\delta,s_0]_s\times \Real_t$.
We first blow up $t=0$, $\rho_1=0$
parabolically (recall that $t$ is the parameter along $\gamma_s$;
$\gamma_s(0)=\mu_1^0(s)$), so that in the interior of
the front face $s$, $t/\rho_1^{1/2}$ and $\rho_1$ become valid coordinates.
Then $(F(\mu_1)-F(\mu_1^0(s)))/\rho_1$ is smooth near the interior of the
front face, and $e^{i(F(\mu_1)-F(\mu_1^0(s)))/\rho_1}$ is smooth on
the blown up space with infinite order vanishing off the front face.
Away from $s=s_0$,
$e^{-is(\sqrt{\mu-\lambda_{1i}-k_2^2/4}-\sqrt{\mu-\mu_1-k_2^2/4})
/\rho_1}$ is rapidly decreasing, hence
the standard push-forward theorem for polyhomogeneous functions
\cite{RBMCalcCon} yields the same result as stationary phase.
In fact, we should think of integrating densities, hence we
need to change $d\mu_1$, i.e.\ $dt$, to $\rho^{1/2}d(t/\rho^{1/2})$, giving
the stationary phase asymptotics $\rho^{1/2}$ times a smooth function
of $(s,\rho^{1/2})$ (multiplied by the exponential $x_1^{iF(\mu_1^0(s))}$
that we took outside the integral in \eqref{eq:integral-77}), though
the terms of the form $\rho$ times a smooth function of $\rho$ and $s$ in fact
cancel (these are the boundary terms of the integration).
This is an alternative way of thinking about complex stationary phase
when the real part of the exponent is non-positive and at least comparable
to or larger than the imaginary part, as suggested to us by Richard Melrose.
It simplifies our previous arguments here, although unfortunately it does
not directly apply when considering asymptotics at the continuous spectrum,
as we do in the next section.

At $s=s_0$ we also need to resolve the geometry of the second factor in
the integrand of \eqref{eq:integral-77}. Had we not performed the parabolic
blow-up above, it would suffice to blow up
$$
F'_1=\{(\rho_1,s,t):\ t=0,\ \rho_1=0,\ s=s_0\}\subset Z_0
$$
in the usual spherical (homogeneous) sense. Indeed, the exponent in the
second factor as well as the denominator are smooth functions of
$(\mu_1-\lambda_{1i})/\rho_1$, hence of $t/\rho_1$ and $(s-s_0)/\rho_1$, and
the exponential is rapidly decreasing as $|t/\rho_1|+|(s-s_0)/\rho_1|\to
\infty$ since
$$
\im(\sqrt{\mu-\lambda_{1i}-k_2^2/4}-\sqrt{\mu-\mu_1-k_2^2/4})
\leq -C'|\mu_1-\lambda_{1i}|\leq -C(|t|+|s-s_0|).
$$
Thus, the second factor is polyhomogeneous on this blown up space,
of order $-1$ on the front face, order $0$ off the front face.

Unfortunately these two blow ups, i.e.\ the parabolic one of $t=0$, $\rho_1=0$,
and the spherical one of $t=0$, $\rho_1=0$, $s=s_0$, conflict with each other,
and we need to find a common resolution. One common resolution is the
following: we blow up the boundary $\rho_1=0$ of $Z_0$ by introducing
$\hat\rho=\rho_1^{1/2}$ as
our new boundary defining function. We write $(Z_0)_{1/2}$ for the
blown up space. The spherical blow-ups of
$$
F_1=\{(\hat\rho,s,t):\ \hat\rho=0,\ t=0\} \Mand
F_2=\{(\hat\rho,s,t):\ \hat\rho=0,\ t=0,\ s=s_0\}
$$
commute with each other since $F_2$ is a p-submanifold of $F_1$, so
$[(Z_0)_{1/2};F_1,F_2]=[(Z_0)_{1/2};F_2,F_1]$.
Let $\ff_j$ denote the front face of
the blow-up of $F_j$, $j=1,2$.
Near the interior of the front face of the second blow-up, i.e.\ $\ff_2$,
$t/\hat\rho$, $(s-s_0)/\hat\rho$ and $\hat\rho$ are valid coordinates.
Now blow up the submanifold
$$
F_3=\{\hat\rho,(s-s_0)/\hat\rho,t/\hat\rho):
\ t/\hat\rho=0,\ (s-s_0)/\hat\rho=0,\ \hat\rho=0\}
$$
(which is a single point) to obtain
$$
Z=[(Z_0)_{1/2};F_1,F_2,F_3].
$$
This introduces, in particular, the front face of the spherical blow up of
$$
F'_1=\{(\rho_1,t,s):\ t=0,\ \rho_1=0,\ s=s_0\}\subset Z_0
$$
in $Z_0$, except that the boundary defining functions differ. In other
words, if we blow up the front face of $[Z_0;F'_1]$ to admit $\rho_1^{1/2}$
as a smooth function, then a neighborhood of the interior of the front face
is naturally diffeomorphic to a neighborhood of the interior of
the front face $\ff_3$ of the blow-up of $F_3$.
This can be seen explicitly since local coordinates which are valid
in this region are given by $(t/\hat\rho)/\hat\rho=t/\rho$,
$((s-s_0)/\hat\rho)/\hat\rho=(s-s_0)/\rho$ and $\hat\rho=\rho^{1/2}$.
Since upon blowing up $F_2$, the lifts of $F_1$ and $F_3$ are disjoint,
their blow-ups commute. Thus, we can rewrite $Z$ as
$[(Z_0)_{1/2};F_2,F_3,F_1]$. Since the first factor of the integrand
of \eqref{eq:integral-77} is polyhomogeneous on $[(Z_0)^{1/2};F_1]$,
while the second factor is polyhomogeneous on $[(Z_0)_{1/2};F_2,F_3]$,
we deduce that the integrand is (one-step) polyhomogeneous on $Z$.

While the geometry of the integrand only requires these blow-ups, we
need further blow-ups to create a b-fibration for the push-forward.
Unfortunately, with the blow-ups discussed above, all we can hope for
is to create a b-fibration with base given by a double blow-up of
$(Y_0)_{1/2}=[0,1)_{\hat\rho}\times(s_0-\delta,s_0]$ (which is
already a blow-up, in the sense of change of the $\Cinf$ structure,
of $Y_0=[0,1)_{\rho_1}\times(s_0-\delta,s_0]$).
Namely, one first blows up
$$
G_2=\{(\hat\rho,s):\ s=s_0,\ \hat\rho=0\},
$$
in $(Y_0)_{1/2}$, and then
$$
G_3=\{(\hat\rho,s):\ (s-s_0)/\hat\rho=0,\ \hat\rho=0\}
$$
in $[(Y_0)_{1/2};G_2]$;
we denote this space by
\begin{equation}\label{eq:Y-def}
Y=[(Y_0)_{1/2};G_2,G_3].
\end{equation}
It is then straight-forward to
carry out appropriate blow-ups in our resolved space $Z$ to obtain a new space
$Z'$ such that $Z'\to Y$ is a b-fibration. Hence we can apply the
push-forward theorem of \cite{RBMCalcCon}, and deduce that the
integral of \eqref{eq:integral-77} yields a polyhomogeneous function
on $Y$. More specifically, the result is a one-step polyhomogeneous
function on $Y$ (keep in mind that $\hat\rho$ is the boundary
defining function!) with order $1$ on the lift of $\hat\rho=0$, and order
$0$ on each of the two new front faces.

The asymptotics from $s\geq s_0$, $s\to s_0$, can be seen similarly.
The optimal contours in $s>s_0$ have been shifted through $\lambda_{1i}$,
so there is automatically a contribution from the pole. We are interested
in letting $s$ decrease to $s_0$. Since in this region the real part of
the exponent
of $x_1^{i\sqrt{\mu_1-k_1^2/4}}$ in the expansion of $R_1(\mu_1)$ is
less than that of $x_1^{i\sqrt{\lambda_{1i}-k_1^2/4}}$, we now need to
expand an operator associated to $H_1$ rather than one associated to $H_2$
around $\lambda_{1i}$. This is a little more delicate, as we discuss below.
The parametrix identity gives
$$
R_1(\mu_1)=P_1(\mu_1)-E(\mu_1)P_1(\mu_1)+E(\mu_1)R_1(\mu_1)F(\mu_1),
$$
with all terms meromorphic. The first two terms are actually holomorphic,
so we can shift the integral to the optimal location, and even let it
go through $\lambda_{1i}$. Thus, we only need to consider the last term.
Here the kernel of $E(\mu_1)$ is polyhomogeneous on the standard double
space $\bM_1\times\bM_1$, with rapid vanishing on the right face.
We expand
the holomorphic function $E_1(\mu_1)$ as
$$
E(\mu_1)=E_1(\lambda_{1i})+(\mu_1-\lambda_{1i})\tilde E_1(\mu_1)
$$
similarly to the expansion for $R_2(\mu-\mu_1)$ above. The contour of
the integral
$\int_{\gamma_s}E_1(\lambda_{1i})R_1(\mu_1)F(\mu_1)R_2(\mu-\mu_1)f\,d\mu_1$
can be shifted farther from $\spec(H_1)$, hence it is negligible compared
to the other terms. So it remains to deal with
$$
\int_{\gamma_s}\tilde E_1(\lambda_{1i})(\mu_1-\lambda_{1i})R_1(\mu_1)
F(\mu_1)R_2(\mu-\mu_1)f\,d\mu_1.
$$
But the integrand above again has the standard asymptotics by virtue
of the holomorphy of $(\mu_1-\lambda_{1i})R_1(\mu_1)$, hence completely
analogous calculation apply as for $s\leq s_0$ above.

We have introduced these new front faces in $Y_0$ in order to understand
the asymptotics in terms of polyhomogeneous expansions. We should
certainly discuss whether this complicated geometry is simply an artifact
of our method, or whether it is necessary in the sense that polyhomogeneous
asymptotics do not hold in a simpler space.
We certainly need at least one of the front faces. Indeed, the residue
term is order $0$ in $s>s_0$, and is missing in $s<s_0$, where the asymptotic
expansion is order $1$ (in terms of $\rho_1^{1/2}$). Thus, it is possible
to obtain polyhomogeneous asymptotics through $s=s_0$ only if there is
some blow-up, and hence a new boundary face, along which the term of order
$0$ tends to zero on approach to the lift of $s<s_0$.  In addition,
the exponential in the second factor is smooth only after $F'_1$
is blown up in $Z_0$, and so it appears that the spherical blow-up of
$s=s_0$, $\rho_1=0$ in $Y_0$, or at least the presence of the
front face of the last blow-up in $[(Y_0)_{1/2};G_2;G_3]$, is necessary.
It then remains to see whether the integral \eqref{eq:integral-77}
is polyhomogeneous on $[Y_0;\{(\rho_1,s):\ \rho_1=0,\ s=s_0\}]$.
If it is, it must be order $0$ on the front face and order $1/2$
on the lift of $s<s_0$. Correspondingly it has order at most $1/2$ on
the blow-up of the corner, which is the interior of the front face
$[(Y_0)_{1/2};G_2]$, although with a different boundary defining function.
Nonetheless, we would expect to see decay in the asymptotics at the
front face of $[(Y_0)_{1/2};G_2]$ away from $G_3$, and there is
no such decay. Indeed, the exponential in the second factor of the
integrand in \eqref{eq:integral-77} is rapidly decreasing in the inverse
image of
this region under the projection, so it can be disregarded. Thus
the dominant term as $\rho_1\to 0$ is an integral of the form
$\int_{\gamma_s} e^{a\mu_1^2/\rho_1}(\mu_1-\lambda_{1i})^{-1}\,d\mu_1$,
$a=-iF''(\mu_1^0(s))/2>0$. Shifting the contour to be vertical,
the integral becomes
$
i\int_{\Real} e^{-at^2/\rho_1}(it+(\mu_1^0(s)-\lambda_{1i}))^{-1}\,dt$,
with imaginary part
$$
\int_{\Real} e^{-at^2/\rho_1}\frac{\mu_1^0(s)-\lambda_{1i}}
{t^2+(\mu_1^0(s)-\lambda_{1i}))^2}\,dt.
$$
Thus, with $\mu_1^0(s)-\lambda_{1i}=\rho_1^{1/2}S$,
$S<0$, and changing variables $T=t/\rho_1^{1/2}$,
$$
\int_{\Real} e^{-aT^2}\frac{S}{T^2+S^2}\,dT<0,
$$
independent of $\rho_1>0$. Hence there is no decay as $\rho_1\to 0$.
The limit $S\to 0$ can also be analyzed; indeed, it is
the distribution $(T+i0)^{-1}$ paired with $e^{-aT^2}$. The limit from
$s>s_0$ is similar, but now we get $(T-i0)^{-1}$ paired with $e^{-aT^2}$.
The difference, $2\pi i$, is exactly the residue corresponding to the pole.

Note moreover that the principal symbol of $R(\mu)f$ in $s< s_0$
is an elliptic $f$-independent
multiple of $P_1^t(\mu_1^0(s))P_2^t(\mu-\mu_1^0(s))f$,
which in turn is the restriction of
$$
x_1^{-k_1/2-i\sqrt{\mu_1^0(s)-k_1^2/4}}R_1(\mu_1^0(s))P_2^t(\mu-\mu_1^0(s))f
$$
to $x_1=0$.
That is, the top term of the asymptotics of $R(\mu)f$,
up to elliptic smooth factors and
after factoring out $x_1^{k_1/2}x_2^{k_2/2}e^{-i\sqrt{\mu-k^2/4}/\rho}$, is
$$
\rho_1^{1/2}x_1^{-k_1/2-i\sqrt{\mu_1^0(s)-k_1^2/4}}R_1(\mu_1^0(s))
P_2^t(\mu-\mu_1^0(s))f,
$$
which behaves as $\frac{\rho_1^{1/2}}{\mu_1^0(s)
-\lambda_{1i}}x_1^{-k_1/2-i\sqrt{\mu_1^0(s)-k_1^2/4}}\Pi_{1i}
P_2^t(\mu-\mu_1^0(s))f$ as $s\to s_0$ (up to lower order terms).
We have thus
explicitly matched the coefficients from the two sides of the corner.

Altogether, we have proved the following theorem.

\begin{thm}\label{thm:res-set-asymp}
Suppose $f\in\dCinf(\bX)$, $\mu\in\Cx\setminus\spec(H)$. Let
$\rho_i=-1/\log x_i$, $\rho^{-1}=\sqrt{\rho_1^{-2}+\rho_2^{-2}}$,
$s=\rho_1/\rho_2$.
Then $R(\mu)f$
has the following asymptotic expansion on $\Xt$:
\begin{equation}\begin{split}
R(\mu)f & =x_1^{k_1/2}x_2^{k_2/2}\exp (-i\sqrt{\mu-k^2/4}/\rho)h \\
 +\sum_{i=1}^{N_1}x_2^{k_2/2+
i\sqrt{\mu-\lambda_{1i}-k_2^2/4}}(\phi_{1i} & \otimes h_{1i})\chi_{1i}
+\sum_{i=1}^{N_2}x_1^{k_1/2+
i\sqrt{\mu-\lambda_{1i}-k_2^2/4}}(h_{2i}\otimes \phi_{2i})\chi_{2i},
\end{split}
\label{eq:asbs}
\end{equation}
where the terms have the following properties.
(See Definition~\ref{Def:Ns-Ps} for the notation, Figure~\ref{fig:contour4}
for a picture.)

\begin{enumerate}
\item
$h_{1i}\in\Cinf(\bM_2)$, $h_{2i}\in\Cinf(\bM_1)$ satisfy
\begin{equation}\begin{split}
&(\phi_{1i}\otimes h_{1i})|_{\bM_1\times \pa\bM_2}
=-(\Pi_{1i}\otimes P_2^t(\mu-\lambda_{1i}))f,\ \text{resp.}\\
&(h_{2i}\otimes \phi_{2i})|_{\pa\bM_1\times \bM_2}.
=-(P_1^t(\mu-\lambda_{2i})\otimes\Pi_{2i})f.
\end{split}\end{equation}

\item
For $\mu$ not real, $\chi_{ji}$ is a function of $s$,
$\chi_{ji}\in\Cinf(\Real^+)$ is identically $1$ or $0$
for sufficiently large and small $s$, with the effect that

\begin{enumerate}
\item
If $\lambda_{1i}\in{\mathcal P}^r_s$, then $\chi_{1i}(s)=0$, and if
$\lambda_{2i}\in{\mathcal P}^\ell_s$, then $\chi_{2i}(s)=0$.

\item
If $\lambda_{1i}\in{\mathcal P}^\ell_s$, respectively if
$\lambda_{2i}\in{\mathcal P}^r_s$, then $\chi_{1i}(s)=1$,
resp.\ $\chi_{2i}(s)=1$.

\item
If $\lambda_{ji}\in{\mathcal N}_s$, then the first term in
(\ref{eq:asbs})
dominates the one corresponding to $\phi_{ji}$, hence the
choice of $\chi_{ji}$ is irrelevant.
\end{enumerate}

\item
If $\mu$ is real, $\mu<\inf\spec(H)$, then $\chi_{1i}$ is a $\Cinf$
function of $(s-s_0)/\rho$, $s_0$ given by
\begin{equation}\label{eq:lambda_1i-def}
\lambda_{1i}=\mu_1^0(s_0)=\frac{\mu-k_2^2/4+s_0^2 k_1^2/4}{1+s_0^2},
\end{equation}
$\chi'_{1i}$ is compactly supported, $\lim_{S\to +\infty}\chi_{1i}(S)=1$,
$\lim_{S\to -\infty}\chi_{1i}(S)=0$. Similarly,
$\chi_{2i}$ is a function
of $(s-s_0)/\rho$, $s_0$ given by
\begin{equation}\label{eq:lambda_2i-def}
\lambda_{2i}=\mu-\mu_1^0(s_0)=\frac{s_0^2(\mu-k_1^2/4)+k_2^2/4}{1+s_0^2},
\end{equation}
$\chi'_{2i}$ is compactly supported, $\lim_{S\to +\infty}\chi_{2i}(S)=0$,
$\lim_{S\to -\infty}\chi_{1i}(S)=1$.

\item
If $\mu\nin\Real$, then $h$ is polyhomogeneous on $\Xt$, of order $1/2$ on the
front face, $3/2$ on the side faces. The principal symbol of $h$ is an
elliptic $f$-independent multiple of
\begin{equation}\begin{split}\label{eq:pr-symbol-88}
&P_1^t(\mu_1^0(s))P_2^t(\mu-\mu_1^0(s))f\ \text{on the front face},\\
&
S_1(k_1^2/4)P_2^t(\mu-k_2^2/4)f\ \text{on the lift of}\ \bM_1\times \pa\bM_2,\\
&P_1^t(\mu-k_1^2/4)S_2(k_2^2/4)f\ \text{on the lift of}\ \pa\bM_1\times \bM_2.
\end{split}\end{equation}

\item
If $\mu\in\Real$, $\mu\nin\spec(H)$, then $g$ is polyhomogeneous on
$\Xt$ doubly blown up at the finite number of submanifolds $s=s_0$, $\rho_1=0$,
of the front face of $\Xt$, as in \eqref{eq:Y-def},
with $s_0$ given by either of the equations
\eqref{eq:lambda_1i-def}-\eqref{eq:lambda_2i-def}.
The order at the front faces of the blow-ups is $0$, while the order
on the old front face is still $1/2$. The principal symbol on the front
faces depends on $f$ only via $(\Pi_{1i}\otimes P_2^t(\mu-\lambda_{1i}))f$,
resp.\ $(P_1^t(\mu-\lambda_{2i})\otimes \Pi_{2i})f$; it is indeed an
elliptic multiple of these.
\end{enumerate}
\end{thm}

\section{Asymptotics inside the continuous spectrum}
The arguments needed to analyze the asymptotics of $R(\mu\pm i0)f$
when $f \in \dCinf(X)$ and $\mu > k^2/4$ are not substantially
different, and in some senses even simpler because the contributions
from the $L^2$ eigenfunctions below the continuous spectrum
are always negligible away from the side-faces. Since threshold eigenvalues
introduce additional complications similar to the discussion of the
previous section, with the additional issue that the push-forward results
of \cite{RBMCalcCon} are not directly applicable since the stationary
phase arguments have a substantially different flavor when the phase
is pure imaginary, we will assume that $k_j^2/4$ is not a pole of $R_j$.
Thus, we have the following theorem.

\begin{thm}\label{thm:cont-spec-asymp}
Suppose $f\in\dCinf(\bX)$, $\mu>k^2/4=(k_1^2+k_2^2)/4$, and $k_j^2/4$ is
not a pole of $R_j$, $j=1,2$. Then, with
$\rho_i=-1/\log x_i$, $\rho^{-1}=\sqrt{\rho_1^{-2}+\rho_2^{-2}}$,
$R(\mu-i0)f$
has the following asymptotic expansion on $\Xt$:
\begin{equation}\begin{split}
& R(\mu-i0)f  =x_1^{k_1/2}x_2^{k_2/2}\exp (-i\sqrt{\mu-k^2/4}/\rho)h \\
& +\sum_{i=1}^{N_1}x_2^{k_2/2+
i\sqrt{\mu-\lambda_{1i}-k_2^2/4}}(\phi_{1i}\otimes h_{1i})
+\sum_{i=1}^{N_2}x_1^{k_1/2+
i\sqrt{\mu-\lambda_{1i}-k_2^2/4}}(h_{2i}\otimes \phi_{2i}),
\end{split}\end{equation}
where $h$ is one-step polyhomogeneous on $\Xt$, and
$h_{1i}$ and $h_{2i}$ are one-step polyhomogeneous
on $\bar M_2$ and $\bar M_1$ respectively. Moreover, $h$ has order
$1/2$ on the front face and $3/2$ on the side faces, while
the $h_{1i}$ is in $\Cinf(\bM_2)$ (i.e.\ has order $0$ at $\pa\bM_2$)
and $h_{2i}\in\Cinf(\bM_1)$. The principal symbol of $h$ is given by
\eqref{eq:pr-symbol-88}.
\end{thm}

\begin{rem}
Note that due to the exponential decay of $x_1^{-k_1/2}
\phi_{1i}$ for non-threshold eigenvalues $\lambda_{1i}$, namely
$$x_1^{-k_1/2}\phi_{1i}
\in x_1^{\sqrt{k_1^2/4-\lambda_{1i}}}\Cinf(\bar M_1),\ k_1^2/4-\lambda_{1i}>0,
$$
and the similar decay of $x_2^{-k_2/2}\phi_{2i}$, the second term is
dominated by the first one everywhere but at the lift of $\bar M_1\times
\pa\bar M_2$, and similarly the third term is dominated
by the first one everywhere but at the lift of $\pa\bar M_1\times
\bar M_2$. In particular, on a neigborhood $U$ of the lift of $\bar M_1\times
\pa\bar M_2$, one has the asymptotics
$$
R(\mu)f|_U=x_1^{k_1/2}x_2^{k_2/2}\exp (-i\sqrt{\mu-k^2/4}/\rho)g+
\sum_{i=1}^{N_1}x_2^{k_2/2+
i\sqrt{\mu-\lambda_{1i}-k_2^2/4}}(\phi_{1i}\otimes g_{1i}),
$$
while on a neighborhood $U'$ of the front face which is disjoint from
$\pa\bar M_1\times \bar M_2$ and $\bar M_1\times
\pa\bar M_2$  both the second and third terms are
irrelevant:
$$
R(\mu)f|_{U'}=x_1^{k_1/2}x_2^{k_2/2}\exp (-i\sqrt{\mu-k^2/4}/\rho)g.
$$
\end{rem}

\begin{rem}
The asymptotics is unchanged away from the corners of $\Xt$ even if
$k_j^2/4$ is a pole of $R_j$, $j=1,2$, in the sense that
these do not contribute to the asymptotics in the interior of the front
face of $\Xt$, and give the standard contribution in each appropriate
side face. This follows from stationary
phase arguments using the integral representation \eqref{eq:tilde-R_2-def}.
In fact, we can always reduce the calculation to the evaluation of an
integral like \eqref{eq:integral-77}, though here it is convenient to
rewrite the fraction appearing in that formula as an integral analogous
to \eqref{eq:tilde-R_2-def}.
\end{rem}

\begin{figure}[ht]
\begin{center}
\mbox{\epsfig{file=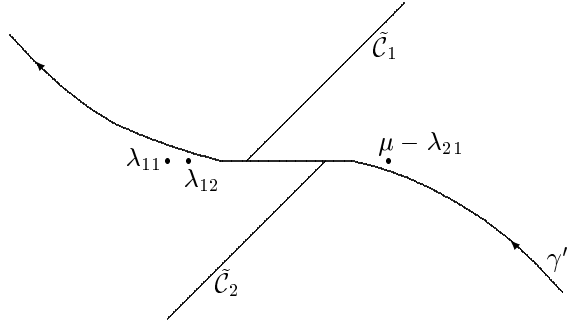}}
\end{center}
\caption{Contour of integration $\gamma'$ to describe $R(\mu)$ inside
the main branch of the continuous spectrum. The cuts ${\mathcal C}_1$ and
${\mathcal C}_2$, corresponding to the spectra of $H_1$ and $H_2$,
respectively, have been rotated to the cuts $\tilde{\mathcal C}_1$ and
$\tilde{\mathcal C}_2$, as in Figure~\ref{fig:contour3}, to make the
picture clearer, although the analytic continuation is not needed
for our proof.}
\label{fig:contour5}
\end{figure}

\begin{proof}
This is simply a limiting case of the argument of the previous section,
but which has additional simplifying features. Thus,
we simply consider a smooth contour $\gamma'$
that runs along the real axis on the interval $[k_1^2/4,\mu-k_2^2/4]$, but
avoids the poles of $R_1$ as well as $\mu$ minus the poles of $R_2$,
see Figure~\ref{fig:contour5}. Such a contour is simply the limit of
the contours $\gamma'$ that we have considered before. Note, in
particular, that with
$$F(\mu_1) =  \sqrt{\mu_1-k_1^2/4}+s\sqrt{\mu-\mu_1-k_2^2/4},
\qquad s=\frac{\rho_1}{\rho_2}$$
as in the previous section,
$\im F(\mu_1)\leq 0$ everywhere, and it is strictly negative when
$\mu_1\nin[k_1^2/4,\mu-k_2^2/4]$. Since the critical point of $F$
lies in $[k_1^2/4,\mu-k_2^2/4]$, the choice of the contour outside this
interval is irrelevant (except that it should be far from the real axis
near infinity, as before, to ensure convergence).
In particular, $\lambda_{ji}\in{\mathcal N}_s$, for all $s\in[0,+\infty)$,
except if $s=0$, when $\im F(\mu-\lambda_{2i})=0$
(but $\im F(\lambda_{1i})<0$); of course a similar statement holds at
$s^{-1}=0$ with $\lambda_{2i}$ replaced by $\lambda_{1i}$. Thus,
the residues of the poles of $R_1$ and $R_2$ also
provide lower order contributions than the stationary phase term in most
regions in $\Xt$, except
that the poles of $R_1$ give the same order as the stationary phase term
at $M_1\times\pa\bM_2$, and similarly for $R_2$.
The standard stationary phase lemma now gives the desired results, much
as in the previous section. To deal with
the points $\mu_1=k_1^2/4$ and $\mu_1=\mu-k_2^2/4$, i.e.\ the endpoints of
$[k_1^2/4,\mu-k_2^2/4]$, one introduces $\tau_1=\sqrt{\mu_1-k_1^2/4}$,
resp.\ $\tau_2=\sqrt{\mu-\mu_1-k_2^2/4}$, just as in the previous section,
and notes that in $0\leq s\leq C$, the phase with respect to $\tau_1$
is never stationary, while with respect to $\tau_2$ it is only stationary
at $s=0$, hence the left end point gives a rapidly decreasing
contribution in this region, while the right end point gives asymptotically
non-trivial terms only at $s=0$, exactly as expected.
\end{proof}

\section{The resolvent compactification}

We now turn to a problem originally posited as one of our main goals,
namely to define a resolution $X^2_{\res}$ of $X^2$, which we call the
resolvent double space, which is a manifold with corners and on which
$R(\mu)$ is polyhomogeneous, at least when $\mu \notin \spec(H)$.
This space is meant to be an analogue of the $0$-double space
$(\bM_j)^2_0$, for either of the conformally compact factors $\bM_j$,
which we discussed in \S 3.  The existence of this space sets the stage
for all further development of the analytic properties of the Laplacian,
including the scattering theory, on $X$.

In fact, the last section contains essentially all of the requisite analysis,
and what remains here is to show how those calculations may be
interpreted. As usual, the starting point is the integral representation
(\ref{eq:resform}) for $R(\mu)$. Again we begin by assuming that neither
operator $H_j$ has bound states.

First, as explained in Theorem~\ref{th:ccst}, apart from the usual diagonal
singularities, the Schwartz kernels of both $R_1(\mu_1)$ and $R_2(\mu-\mu_1)$
are polyhomogeneous at the boundary hypersurfaces of the $0$-double spaces
$(\bM_j)_0^2$. Hence the Schwartz kernel of the integrand $R_1(\mu_1)
R_2(\mu-\mu_1)$ is polyhomogeneous on the product of these $0$-double spaces,
$(\bM_1)^2_0 \times(\bM_2)^2_0$, again with diagonal singularities on each
factor, and thus it is clear that our resolution process should start here.
Our task is to see what new singularities are introduced by the
contour integration, and how these may be accomodated geometrically.
We shall systematically neglect the diagonal singularities in the ensuing
discussion; that they produce the correct singularity upon integration is
obvious in the interior, and requires only minor justification, which we omit.

This product of double spaces has six boundary hypersurfaces, namely
$$
\ff_1\times (\bM_2)_0^2,\ \lf_1\times (\bM_2)_0^2,\ \rf_1\times (\bM_2)_0^2,
\ (\bM_1)_0^2\times\ff_2,\ (\bM_1)_0^2\times\lf_2,\ (\bM_1)_0^2\times\rf_2,
$$
where $\ff_j$ is the front face, and $\lf_j$ and $\rf_j$ the left and right
faces of $(\bM_j)^2_0$. The kernel of $R_j(\nu)$ is the product of a
smooth function on $(\bM_j)^2_0$ (albeit with a diagonal singularity)
with the factor
$$
(\rho_{\lf_j}\rho_{\rf_j})^{k_j/2+i\sqrt{\nu-k_j^2/4}};
$$
in particular, it is smooth up the $\ff_j$. It is straightforward to see
that $R(\mu)$ is smooth at the two `front faces' of the product,
$\ff_1 \times (\bM_2)_0^2$ and $(\bM_1)_0^2 \times \ff_2$. This is because
the exponents in the expansions of $R_1(\mu_1) R_2(\mu-\mu_1)$ at these
faces do not depend on $\mu_1$. It remains to analyze the behaviour
at the remaining `side faces'.

In fact, the computations of the last section give the asymptotics of
$R(\mu)$ at $(\lf_1\times (\bM_2)_0^2)\cup((\bM_1)_0^2\times\lf_2)$, and
in particular at the corner $(\lf_1\times(\bM_2)_0^2)\cap((\bM_1)_0^2
\times\lf_2)$. To proceed further, the salient observation is that while
those stationary phase computations were motivated by the geometric
picture of letting the variable $(z_1,z_2) \in M_1 \times M_2$ converge to
the boundary (in one of three ways), they can also be interpreted purely
analytically as a derivation of asymptotics when integrating a function
of the form $x_1^{\alpha_1}x_2^{\alpha_2}F$, where $\alpha_1$ and $\alpha_2$
are the appropriate exponents depending on $\mu$ and $\mu_1$ and $F$ depends
smoothly on these variables. What we are asking now is the asymptotics of
an integral of almost exactly the same form, but where each $x_j$ is
replaced by $\rho_{\lf_j}\rho_{\rf_j}$.

Replace each of the exponential
type defining function $\rho$ at the side faces of these $0$-double spaces
by the polynomial type defining function $-1/\log \rho$. In other words,
logarithmically blow up each of the two side faces $\lf_j$ and $\rf_j$ of
$(\bM_j)_0^2$, but {\it not} the front face. Call the resulting space
$(\bM_j)^2_{0,\log}$. For notational convenience, set
\[
R_{\lf_j} = -1/\log \rho_{\lf_j},\qquad R_{\rf_j} = -1/\log \rho_{\rf_j},
\qquad j = 1,2.
\]
Now perform the same calculations as in the last section. It follows
that the asymptotics of $R(\mu)$ along the interiors of each of the four
side faces is polyhomogeneous in terms of these new defining functions.

However, to fully resolve the singularities at the various corners, we must
perform a further sequence of blow-ups. More specifically, the stationary
phase calculation produces terms polyhomogeneous in $R_{\lf_j}$ and
$R_{\rf_j}$, and the only other terms which arise here are polyhomogeneous
in $(\log \rho_{\lf_2}\rho_{\rf_2})/(\log \rho_{\lf_1}\rho_{\rf_1})$.
Thus we need to find a resolution so that this projective coordinate
is $\mathcal C^\infty$ when it is bounded, and similarly for its inverse.
It is not hard to do this. We can certainly assume that $\rho_{\lf_j},
\rho_{\rf_j} \leq 1/2$. Now observe that
\begin{equation}\label{eq:log-rho-id-9}
\frac{\log(\rho_{\lf_2}\rho_{\rf_2})}{\log (\rho_{\lf_1}\rho_{\rf_1})}=
\frac{\log\rho_{\lf_2}}{\log \rho_{\lf_1}+\log\rho_{\rf_1}}
+\frac{\log\rho_{\rf_2}}{\log \rho_{\lf_1}+\log\rho_{\rf_1}},
\end{equation}
which is equal to
\begin{equation}\label{eq:log-rho-id-99}
\frac{R_{\lf_1}}{R_{\lf_2}} \frac{R_{\rf_1}}{R_{\lf_1}+R_{\rf_1}}
+ \frac{R_{\lf_1}}{R_{\rf_2}} \frac{R_{\rf_1}}{R_{\lf_1}+R_{\rf_1}}.
\end{equation}
Both terms in this last expression are positive, and hence their sum
is bounded if and only if both terms are. Clearly we require a space
where all of the quotients
\[
R_{\lf_i}/R_{\lf_j}, \quad R_{\lf_i}/R_{\rf_j}, \quad R_{\rf_i}/R_{\rf_j}
\qquad i,j = 1,2,
\]
and their inverses are smooth (when bounded).

Thus first blow up
\begin{equation}\label{eq:double-blow-up-1}
(\lf_1\cap\rf_1)\times (\bM_2)^2_{0,\log},\ (\bM_1)^2_{0,\log}
\times(\lf_2\cap\rf_2).
\end{equation}
This may be done in either order since these submanifolds are transverse.
Next, we must blow up the collection of submanifolds covering
\begin{equation*}\begin{split}
&(\lf_1\times (\bM_2)^2_{0,\log})\cap((\bM_1)^2_{0,\log}\times\lf_2),
\ (\lf_1\times (\bM_2)^2_{0,\log})\cap((\bM_1)^2_{0,\log}\times\rf_2),\\
&(\rf_1\times (\bM_2)^2_{0,\log})\cap((\bM_1)^2_{0,\log}\times\lf_2),
\ (\rf_1\times (\bM_2)^2_{0,\log})\cap((\bM_1)^2_{0,\log}\times\rf_2).
\end{split}\end{equation*}
These are disjoint now after the initial blow-up, and so may be blown up in
any order. The manifold $X^2_{\res}$ obtained in this way has all the desired
properties. In fact, working in a region where \eqref{eq:log-rho-id-9} is
bounded, we must show that both terms in \eqref{eq:log-rho-id-99} are smooth.
But
\[
\frac{R_{\lf_1}}{R_{\lf_2}} \frac{R_{\rf_1}}{R_{\lf_1}+R_{\rf_1}}
=\frac{R_{\lf_1}}{R_{\lf_2}} \frac{1}{(R_{\lf_1}/R_{\rf_1}) + 1},
\]
with a similar expression for the other term. The first spherical blow up
\eqref{eq:double-blow-up-1} makes $R_{\rf_j}/R_{\lf_j}$ and its inverse
smooth where bounded, while the prefactors $R_{\lf_1}/R_{\lf_2}$ and
$R_{\lf_1}/R_{\rf_2}$ also become smooth because $(\lf_1 \times
(\bM_2)^2_{0,\log}) \cap ((\bM_1)^2_{0,\log} \times \lf_2)$ and
$(\lf_1 \times (\bM_2)^2_{0,\log}) \cap ((\bM_1)^2_{0,\log} \times \rf_2)$
have been blown up. We thus conclude that when $H_1$ and $H_2$ have no
$L^2$ eigenfunctions, $R(\mu)$ is polyhomogeneous on this space
$X^2_\res$.

To conclude this discussion, suppose that either $H_1$ or $H_2$ have
bound states. Then $R(\mu)$ will have additional
singularities at certain submanifolds of $X^2_\res$ where
$\log(\rho_{\lf_2}\rho_{\rf_2})/\log (\rho_{\lf_1}\rho_{\rf_1})$ assumes
the same critical values $s_0$ as already appeared for $s=\log x_2/\log x_1$
in part (v) of Theorem~\ref{thm:res-set-asymp}. That is, these singularities
happen at the submanifolds given by $S=s_0$ {\em inside $\pa X^2_{\res}$},
where $S$ is the smooth function
on $X^2_\res$ given by
$$
S=\left(\frac{R_{\lf_1}}{R_{\lf_2}}+ \frac{R_{\lf_1}}{R_{\rf_2}}\right)
\frac{1}{(R_{\lf_1}/R_{\rf_1}) + 1}=(\sigma_{l1l2}+\sigma_{l1r2})
\frac{1}{\sigma_{l1r1}+1},
$$
in the region where $\sigma_{l1l2}=\frac {R_{\lf_1}}{R_{\lf_2}}$,
$\sigma_{l1r2}=\frac{R_{\lf_1}}{R_{\rf_2}}$
and $\sigma_{l1r1}=R_{\lf_1}/R_{\rf_1}$ are valid coordinates,
i.e.\ where they
are bounded. Since the second factor is thus bounded from both below and from
above, we conclude that, as long as $s_0=0$ is not one of the critical values,
at least one of $\sigma_{l1l2}$ and $\sigma_{l1r2}$ is non-zero at any point
on $S=s_0$, and in addition
$\pa_{\sigma_{l1l2}}S\neq 0$, $\pa_{\sigma_{l1r2}}S\neq 0$ there, so $S=s_0$
is a codimension one p-submanifold of $X^2_{\res}$, and we have
singularities at its intersection with $\pa X^2_{\res}$, which
in turn is
a finite union of codimension one p-submanifolds of $\pa X^2_{\res}$.
These singularities can be resolved by the same double blow up, preceeded
by the square root blow up of the boundary defining function, as
on the single space.

While these blow ups may appear somewhat complicated, the only relevant
part of $X^2_{\res}$ as far as the Martin boundary is concerned,
is the inverse image $U'$ of regions $X\times U$, $\overline{U}$ compact
subset of $X$, $U=U_1\times U_2\subset \bM_1\times \bM_2$ open,
under the blow-down map $\beta:X^2_{\res}\to X^2$.
Thus, in the second factor we always stay away from the boundary,
hence any blow-ups involving $\rf_1$, $\rf_2$ can be neglected.
Thus, the intersection of $U'$ with $X^2_{\res}$ is a subset of
$$
[(\bM_1)_{\log}\times U_1\times (\bM_2)_{\log}\times U_2;
\pa (\bM_1)_{\log}\times U_1\times\pa (\bM_2)_{\log}\times U_2]
$$
which is essentially the same as $\Xt\times U_1\times U_2$ (i.e.\ they agree
up to the rearranging of factors). Hence, in this set the asymptotics
of the kernel of $R(\mu)$ is
given by Theorem~\ref{thm:res-set-asymp}, with an extra variable on $U$
added, and references to $f$ dropped (this really corresponds to applying
the resolvent to delta distributions at points $z$ and letting $z$ vary).

\section{The Martin compactification}

With the information we have collected and proved about the resolvent
$R(\mu)$, it is now a simple matter to determine the full Martin
compactification of $X = M_1 \times M_2$.
We refer back to \S 5 for the general details of how this construction is
to be carried out, but we briefly recall the main ideas.
Fix any $\mu \in \RR$, $\mu<\inf\spec(H)$. In addition,
let $p\in X$ be fixed and $w^{(\ell)} = (\well_1,\well_2) \in X$ any sequence
which leaves every compact set. Define
\[
U^{(\ell)}(z) = \frac{R(\mu;z,\well)}{R(\mu;p,\well)};
\]
this function is a solution of $(H - \mu)U^\ell = 0$ for $z \neq \well$
and is normalized so that $U^{(\ell)}(p) = 1$. Any such sequence
has a subsequence $U^{\ell'}$ which converges uniformly on compact
sets to a function $U$ which satisfies $(H-\mu)U = 0$ on all of $X$
and $U(p) = 1$. The Martin boundary $\overline{X}_M$ is defined to be the
set of all possible limit eigenfunctions that arise in this manner.

First consider the case if neither $H_1$ nor $H_2$ has $L^2$ eigenvalues.
In this case, by the asymptotics of Proposition~\ref{prop:no-bd-state-asymp}
(see also the last paragraph of the previous section)
\begin{equation}\begin{split}
R(\mu)|_{\Xt\times X}&=x_1^{k_1/2}x_2^{k_2/2}
\exp (-i\sqrt{\mu-k^2/4}/\rho)h_0 g,
\end{split}\end{equation}
where $h_0$ is the product of the square root of a defining function
of the front face of $\Xt$ with the $3/2$ power of a defining function
of each side face of $\Xt$, and $g$ is $\Cinf$ in $\Xt\times X$ with
restriction to each face being an elliptic multiple of
\begin{equation}\begin{split}
&P_1^t(\mu_1^0(s))P_2^t(\mu-\mu_1^0(s))\ \text{on the front face},\\
&
S_1(k_1^2/4)P_2^t(\mu-k_2^2/4)\ \text{on the lift of}\ \bM_1\times
\pa\bM_2\times U,\\
&P_1^t(\mu-k_1^2/4)S_2(k_2^2/4)\ \text{on the lift of}\ \pa\bM_1\times \bM_2
\times U.
\end{split}\end{equation}
In particular,
$$
U(\mu,z,w)=\frac{R(\mu;z,w)}{R(\mu;p,w)}
$$
is a continuous function on $\Xt_w\times X_z$.
Thus, the map
\begin{equation}\label{eq:Xt-to-XM}
\pa\Xt\ni w\mapsto U(\mu,.,w)\in\Cinf(X)
\end{equation}
defines a map $\Xt\to\pa\overline{X}_M$ which is continuous and surjective.
More explicitly, if $\well$ is any sequence of points in $X$ converging
to $w\in\pa\Xt$, the continuity of $U$ implies that
$U^{(\ell)}(z)=U(\mu,z,w)$, hence $U(\mu,.,w)$ is a point in
the Martin boundary. Conversely,
suppose that $\well$ is any sequence of points in $X$
which leaves every compact subset of $X$, and such that $U^{(\ell)}$
converges uniformly on compact subsets of $X$. By passing to a subsequence
$w^{(\ell')}$ we may assume that $w^{(\ell')}$ converges to some $w\in\pa
\Xt$ since $\Xt$ is compact. Hence $\lim U^{(\ell')}(z)=U(\mu,w,z)$,
so the map is indeed surjective.

Explicitly, $U$ is given by
\begin{equation}\begin{split}
&\frac{P_1^t(\mu_1^0(s);z_1,w_1)P_2^t(\mu-\mu_1^0(s);z_2,w_2)}
{P_1^t(\mu_1^0(s);p_1,w_1)P_2^t(\mu-\mu_1^0(s);p_2,w_2)}
,\ (s,w_1,w_2)\in
[0,+\infty)_s\times\pa \bM_1\times\pa\bM_2,\\
&
\frac{S_1(k_1^2/4;z_1,w_1)P_2^t(\mu-k_2^2/4;z_2,w_2)}
{S_1(k_1^2/4;p_1,w_1)P_2^t(\mu-k_2^2/4;p_2,w_2)},
\ (w_1,w_2)\in\bM_1\times \pa\bM_2,\\
&\frac{P_1^t(\mu-k_1^2/4;z_1,w_1)S_2(k_2^2/4;z_2,w_2)}
{P_1^t(\mu-k_1^2/4;p_1,w_1)S_2(k_2^2/4;p_2,w_2)},
\ (w_1,w_2)\in\pa\bM_1\times \bM_2.
\end{split}\end{equation}
Thus, the injectivity of the map \eqref{eq:Xt-to-XM} depends on whether
these are all different elements of $\Cinf(X_z)$ as $w$ varies.
It is easy to see that the restriction of this map to the front face is
indeed injective since $P_j^t(\nu;z_j,w_j)$ has a pole when
$z_j \to w_j$ and is otherwise bounded, and so the points $w_1$ and
$w_2$ are uniquely identified, as is the eigenparameter $\mu_1^0(s)$
and hence $s$ itself. On the other hand, it is not obviously true that
the generalized spherical function $S_j(k_j^2/4;z_j,w_j^0)$ determines
the point $w_j^0$. In case $M_j$ is hyperbolic space, this is holds
because of the rotational symmetry around $w_j^0$, and so must also
hold for conformally compact metrics (globally) near to the hyperbolic one.
In general, these portions of the Martin boundary are identified with
${\mathcal I}_1 \times \bM_2$ and $\bM_1 \times {\mathcal I}_2$, where
\[
\mathcal I_j = \mbox{image}\big( w_j^0 \longrightarrow
S_j(k_j^2/4;z_j,w_j^0)\big).
\]

Next suppose that $H_1$ has $L^2$ eigenvalues.
The asymptotics of the resolvent kernel in $\Xt\times X$ are now
given by Theorem~\ref{thm:res-set-asymp}, again in the sense discussed
in the last paragraph of the previous section. Let $\lambda_{11}$ be
the bottom eigenvalue of $H_1$ which is hence simple, and let
$s_0$ be the corresponding value of $s=\rho_1/\rho_2$ given by
\eqref{eq:lambda_1i-def}. In the region $s<s_0$, in particular near
$\pa \bM_1\times\bM_2\times X$, $R(\mu)$ has the same asymptotics as
before, while on the front face of the double blow up of $s=s_0$ as
well as in $s>s_0$ the leading part
of the asymptotics is given by elliptic multiples of
$x_1^{-k_1/2-i\sqrt{\lambda_{1i}-k_1^2/4}}
\Pi_{11}\times P_2^t(\mu-\lambda_{1i})$ which are independent of $z\in X$.
We again set
$$
U(\mu,z,w)=\frac{R(\mu;z,w)}{R(\mu;p,w)},
$$
and $U$ is again continuous on $\Xt\times X$. Now, however, $w\mapsto
U(\mu,.,w)$ factors through the map $\pa\Xt\to \pa\Xt_c$, where $\pa\Xt_c$
is the collapsed boundary of $\Xt$ given by
$$
\pa\Xt_c=(\pa \bM_1\times\bM_2)\cup ([0,s_0)_s\times \pa\bM_1\times\pa\bM_2)
\cup \pa\bM_2,
$$
and the collapse map $\pa\Xt\to\pa\Xt_c$ is the identity on the first
two sets on the right hand side and is the projection $[s_0,+\infty)
\times\pa\bM_1\times\pa\bM_2\to\pa\bM_2$, resp.\ $\bM_1\times\pa\bM_2
\to\pa\bM_2$ in the other parts of $\pa\Xt$.
Thus, we obtain a continuous surjection from $\pa\Xt_c$ to
$\pa\overline{X}_M$.

Again, $U$ can be written down explicitly. It is the same as in the
case without $L^2$ eigenfunctions when $s<s_0$, and it is
$$
\frac{P_2^t(\mu-k_2^2/4;z_2,w_2)}
{P_2^t(\mu-k_2^2/4;p_2,w_2)},
\ w_2\in \pa\bM_2.
$$
Thus, the map $\pa\Xt_c\to\pa\overline{X}_M$ is injective over the
collapsed part of the boundary, as can be seen by letting
$z_2\to w_2$, and otherwise we are exactly in the same situation as
in the case where there are no $L^2$ eigenfunctions.

Similar arguments work if $H_2$ has eigenvalues, or even if both $H_1$ and
$H_2$ do. In the latter case $\Xt$ collapses near both side faces.

\bibliographystyle{plain}
\bibliography{sm}

\end{document}